\documentclass[twocolumn]{autart}    

\usepackage{graphicx}          
\usepackage{amssymb}
\usepackage{aligned-overset}
\usepackage{amsmath} 
\usepackage{subfigure}

\usepackage{algorithm}
\usepackage{etoolbox}

\let\classAND\AND
\let\AND\relax
\usepackage{algorithmic}

\let\AND\classAND
\AtBeginEnvironment{algorithmic}{\let\AND\algoAND}
\usepackage{color}

\newtheorem{problem}{Problem}
\newtheorem{lemma}{Lemma}

\begin{document}

\begin{frontmatter}

\title{Nonlinear Optimal Guidance with Constraints on Impact Time and Impact Angle\thanksref{footnoteinfo}} 

\thanks[footnoteinfo]{This paper was not presented at any IFAC 
meeting. Corresponding author Zheng Chen. Tel. +86-13149205555.}


\author[1]{Fanchen Wu}\ead{fanchen-w@zju.edu.cn},    
\author[1,2,3]{Zheng Chen}\ead{z-chen@zju.edu.cn},               
\author[1,3]{Xueming Shao}\ead{mecsxm@zju.edu.cn},  
\author[1]{Kun Wang}\ead{wongquinn@zju.edu.cn}

\address[1]{School of Aeronautics and Astronautics, Zhejiang University, Hangzhou 310027, Zhejiang, China}  
\address[2]{Huanjiang Laboratory, Zhuji 311899, Zhejiang, China}           
\address[3]{State Key Laboratory of Fluid Power Mechatronic Systems, Hangzhou 310027, China}        

\begin{keyword}                           
Nonlinear optimal guidance; Artificial neural network; Impact time; Impact angle.               
\end{keyword}                             

\begin{abstract}                         
This paper aims to address the nonlinear optimal guidance problem with impact-time and impact-angle constraints, which is fundamentally important for multiple pursuers to collaboratively achieve a target. Addressing such a guidance problem is equivalent to solving a nonlinear minimum-effort control problem in real time. To this end, the Pontryagain's maximum principle is employed to convert extremal trajectories as the solutions of a parameterized differential system. The geometric property for the solution of the parameterized system is analyzed, leading to an additional optimality condition. By incorporating this optimality condition and the usual disconjugacy condition into the parameterized system, the dataset for optimal trajectories can be generated by propagating the parameterized system without using any optimization methods. In addition, a scaling invariance property is found for the solutions of the parameterized system. As a consequence of this scaling invariance property, a simple feedforward neural network trained by the solution of the parameterized system, selected at any fixed time, can be used to generate the nonlinear optimal guidance within milliseconds. Finally, numerical examples are presented, showing that the nonlinear optimal guidance command generated by the trained network can not only ensure the expected  impact angle and impact time are precisely met but also requires less control effort compared with existing guidance methods.
\end{abstract}

\end{frontmatter}

\section{Introduction}
The capability of a pursuer to intercept a target is affected by its guidance law as it determines trajectory or control command.
For this reason, various guidance laws have been developed in the literature since the 1960s. 
Traditional guidance laws are usually designed to ensure the terminal miss distance is as close as possible to zero. In recent decades, additional terminal constraints have been incorporated into guidance laws in order for better effect of interception. For instance, constraints on impact time and impact angle are important for multiple pursuers to cooperatively achieve a target \cite{2006Insoo}. To this end, Time- and Angle-Constrained Guidance (TACG) laws have been extensively studied in the literature.
In general, methods for designing TACG laws can be categorized into three groups: 1) Proportional Navigation (PN) based methods, 2) advanced control theory based methods, and 3) trajectory shaping based methods. 

As the PN, probably the most popular guidance law, is well known to be simple and robust to implement, researchers naturally designed TACG laws by introducing additional terms into the conventional PN. Jeon {\it et al}. in \cite{2006Insoo} developed an Impact-Time-Control Guidance (ITCG) law by introducing the difference between estimated time-to-go and expected time-to-go into the PN.
Jeong {\it et al}. in \cite{2004Kim} developed an Impact-Angle-Control Guidance (IACG) law, by including a bias term associated with the Line-of-Sight (LOS) rate error, to control the impact angle. 
To simultaneously control both impact time and impact angle, Lee {\it et al}. in \cite{2007Insoo} introduced a bias term for impact time errors in a minimum-jerk guidance law to precisely adjust the impact time and impact angle. Considering that jerk served as the control command and the lateral acceleration command was derived via time integration of this command in \cite{2007Insoo},
Zhang {\it et al}. in \cite{2013Zhang} developed a simplified form of TACG that directly employs the pursuer's normal acceleration as the control command. 
Recently, in addition to the impact-time and impact-angle constraints, the field-of-view limitation was taken into account in \cite{2022Chenye} by incorporating boundary conditions of the field-of-view into the control command.
Note that the biased terms related to impact angle in the above guidance laws \cite{2007Insoo,2013Zhang,2022Chenye} were derived under the assumption that the collision course is close to a nominal one. Therefore, those guidance laws cannot preserve optimality once the deviations from the collision triangle are relatively large \cite{2019ChenZ}.

Apart from using PN-based laws to control impact time and impact angle, advanced control theories, such as Lyapunov theory \cite{2016Saleem,KIM20142509} and Sliding Mode Control (SMC) \cite{2012Harl,2020Han,2016Zhao,2019Chenxiaot} have also been employed. Saleem and Ratnoo \cite{2016Saleem} utilized the Lyapunov theory to design an ITCG law so that the terminal impact time can be controlled by tuning a single parameter in an exact expression of impact time. 
Similarly, Kim {\it et al}. in \cite{KIM20142509} proposed to control impact angle by augmenting impact-angle error to a Lyapunov candidate function. Nevertheless, it is rare to see studies on using Lyapunov theory to control both impact time and impact angle. 

As to the SMC-based guidance laws, the usual way is to first design a sliding surface. Then, the guidance law is derived by considering driving the pursuer to slide along the surface.
For instance, an LOS rate shaping process was introduced by Harl {\it et al}. \cite{2012Harl} to create such a sliding surface, and an SMC law was developed to track the desired line-of-sight rate profile. This method was later extended to the 3D scenario by Han {\it et al}. \cite{2020Han}. A specific time-varying sliding surface with two unknown coefficients was designed by Zhao {\it et al}. \cite{2016Zhao}. These two coefficients need to be independently tuned to meet the constraints of impact time and impact angle.
The coefficients in \cite{2012Harl,2016Zhao} often require trial-and-error tuning or optimization routines, leading to substantial computational efforts.
To avoid the need for parameter tuning, Chen {\it et al}. \cite{2019Chenxiaot} introduced a novel SMC-based  guidance law. 
To be specific, the pursuer was initially guided to a sliding surface within a finite time and maintained its position on the surface until the impact time.
To summarize, while the guidance laws based on SMC can be applied to address problems with diverse uncertainties and complex constraints, their performance can be significantly influenced by the choice of the sliding surface and the tuning of guidance parameters. Additionally, the chattering phenomenon is difficult to eliminate in most of the SMC-based guidance laws \cite{2006Vadim}.

Another popular approach is to use trajectory shaping to control impact time and impact angle. The control actions are typically formulated as high-order polynomials, and the coefficients of these polynomials are determined to meet the complex boundary conditions of the engagement.
For example, the guidance command was defined as a polynomial function with three unknown coefficients in \cite{2013Kim}. While two of the coefficients were determined to satisfy the impact angle, the remaining one was adjusted to control the impact time. For simplicity, the coefficient related to impact time was obtained using bang-bang control.
Similarly, Tekin {\it et al}. in \cite{2016Tekin} designed a three-phased open-loop control profile for simultaneous control of impact time and impact angle, and a brute-force numerical search method was used to adjust the parameters of the profile.
In addition to directly shaping the control command, Kang {\it et al}. in \cite{2019Kang} addressed the problem by shaping the look angle as a generalized polynomial, and close-to-optimal solutions were attained for low-order polynomials.

In order to simultaneously control impact time and impact angle, multi-stage guidance laws have also been developed in the literature (see, e.g., \cite{2012Harrison,2015Kumar2,2017Song,2018Hu}) by switching between trajectory-shaping guidance laws, SMC-based guidance laws and some other approaches.
Whereas, the performance of such guidance laws is significantly affected by the switching logic, and it is not easy to select an optimal switching logic. To ensure that a pursuer moves in an optimal way while strictly satisfying impact-time and impact-angle constraints, one needs to design optimal guidance laws that should not only satisfy complex constraints (e.g.,  impact-time and impact-angle constraints) and the nonlinear kinematics but also optimize a specific performance index. For this reason, the optimal guidance problem with nonlinear kinematics, dubbed as Nonlinear Optimal Guidance (NOG), has become an active research topic in the last decade. 

In order to address an NOG problem, one usually needs to develop a real-time method for a corresponding nonlinear optimal control problem. However, it is well known that existing numerical methods cannot be guaranteed to compute the optimal solution in real time as they usually suffer the issue of convergence. In order to get real-time solutions, the nonlinear kinematics is usually simplified or linearized to get closed-form solutions. For example, 
the classic linear kinematic approximation was used for a series of optimal guidance problems in \cite{1998Shneydor} (Chapter 8). Merkulov {\it et al.} in \cite{2022Merkulov} employed the quadratic approximation on equations of motion to address the minimum-effort guidance problem with fixed-terminal time.
The first attempt for a closed-form solution for NOG was presented in \cite{1984Guelman}.
In this work, a set of multiple nonlinear equations was solved numerically.
However, numerical solvers cannot find all roots of multiple nonlinear equations, potentially resulting in a local optimum instead of the desired global optimum.
To find the global optimal solution for the NOG described in \cite{1984Guelman}, a parameterization method was developed by Chen and Shima in \cite{2019ChenZ}. This approach involved generating semi-analytical solutions by identifying the zeros of a real-valued function.
As a natural extension of \cite{2019ChenZ}, Wang {\it et al.} developed a parameterized approach for addressing the NOG problem with fixed impact time \cite{2022Wang}. This parameterized approach was directly applied in \cite{Cheng2023} to the NOG problem engaged in the vertical plane. However, to the best of the authors' knowledge, research on NOG with constraints on both impact time and impact angle is scarce to see in the literature. 

For simplicity of presentation, we shall use TAC-NOG to denote the Time- and Angle-Constrained NOG.
The TAC-NOG command is actually determined as the solution of a nonlinear minimum-effort control problem with fixed terminal impact time and fixed terminal impact angle. If one is able to solve the corresponding minimum-effort control problem in real time, the TAC-NOG command can be obtained immediately without using time-to-go estimation or parameter tuning \cite{2022Wang}. 
To achieve this, a parameterized system is established in the current paper by using the Pontryagin's Maximum Principle (PMP).
In addition, some additional optimality conditions and geometric properties are derived for the optimal trajectories related to the TAC-NOG command. By incorporating the optimality conditions and geometric properties into the parameterized system,
one is allowed to construct the dataset for the mapping from state to optimal feedback control command via simply propagating some differential equations. Furthermore, a scaling invariance property is found for the solutions of the parameterized system, enabling an Artificial Neural Network (ANN) trained by the dataset, selected at any fixed time, to generate the TAC-NOG commands in real time.

The remaining sections of this paper are structured as follows. The nonlinear minimum-effort control problem is formulated in Section \ref{Se:formulation}. The procedure for generating real-time solutions via ANN is presented in Section \ref{Se:Real-Time}. In Section \ref{Se:Characterizations}, some optimality conditions are derived, and a parameterized system is formulated for constructing dataset. Numerical examples are presented in Section \ref{Se:Simulations}.

\section{Preliminary}\label{Se:formulation}
\subsection{Problem Formulation}

\begin{figure}
\begin{center}
\includegraphics[height=4cm]{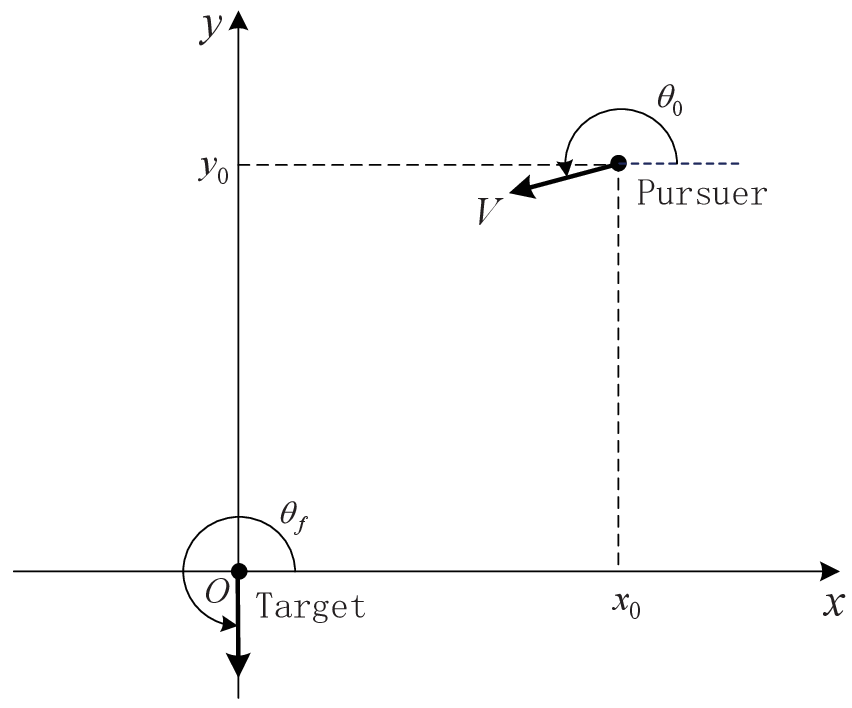}    
\caption{Two-dimensional coordinate system.}  
\label{Fig:geometry}                                 
\end{center}                                 
\end{figure}

Consider a pursuer and a static target on a 2-dimensional coordinate system Oxy, as shown in Fig.~\ref{Fig:geometry}. 
The frame Oxy has its origin located at the target, with the positive x-axis oriented towards the East and the positive y-axis pointing Northward.
Denote by $\boldsymbol{z}=(x,y,\theta)$ the state of the pursuer, comprising a position vector $(x,y)\in \mathbb{R}^2$ and a heading angle $\theta\in [0,2\pi)$ which is positive when measured counterclockwise. The speed of the pursuer, denoted by $V$, is considered to be constant. Then, the nonlinear kinematics of the pursuer can be expressed as \cite{2006Luping}
\begin{equation}
    \label{Eq:system1_2}
    \left\{
        \begin{aligned}
		\dot{x} = & V\cos \theta(t) \\
		\dot{y} = & V\sin \theta(t) \\
		\dot{\theta} = & \dfrac{u(t)}{V}
		\end{aligned}
	\right.
\end{equation}
where $t\in \mathbb{R}_+$ denotes time, the over dot denotes the differentiation with respect to time, and $u$ is the control command.
By normalizing the speed to one, the kinematics can be simplified to
\begin{equation}
    \label{Eq:system1}
    \left\{
        \begin{aligned}
		\dot{x} = & \cos \theta(t) \\
		\dot{y} = & \sin \theta(t) \\
		\dot{\theta} = & u(t)
		\end{aligned}
	\right.
\end{equation}

Then, the real-time generation of the TAC-NOG command is equivalent to addressing the following Minimum-Effort Control Problem (MECP) in real time.
\begin{problem}\label{Problem1}
Given an initial condition $\boldsymbol{z}(0)=(x_0,y_0,\theta_0)$, a fixed terminal time $t_f\in \mathbb{R}_+$, and a terminal heading angle $\theta_f \in [0,2\pi)$, the MECP consists of steering the system in Eq. (\ref{Eq:system1}) with a measurable control $u(\cdot)$ on the interval $[0,t_f]$ from the initial state $\boldsymbol{z}(0)$ to the final state $\boldsymbol{z}(t_f)=(0,0, \theta_f ) $ while minimizing the control effort, i.e.,
\begin{equation}
min \ J=\int_{0}^{t_f} \dfrac{1}{2}u^2\, dt
\end{equation}
\end{problem}

Without loss of generality, we assume the terminal
heading angle in Problem \ref{Problem1} to be fixed as $-\pi/2$, i.e., $\theta_f = -\pi/2$. It's worth noting that if the final heading angle differs from $-\pi/2$, one can use a coordinate rotation to change the final heading angle as $-\pi/2$.
The procedure for rotating the trajectory governed by Eq. (\ref{Eq:system1}) is given in Appendix \ref{App:A}.

\subsection{Necessary Conditions} 
Denote by $\boldsymbol{p}=(p_x,p_y,p_\theta) \in \mathbb{R}^3$ the costate variables of $\boldsymbol{z}=(x,y,\theta)$. Then, the Hamiltonian of the MECP is expressed as
\begin{align}
	H = p_x \cos\theta + p_y \sin\theta + p_\theta u - \dfrac{1}{2} u^2
\end{align}
The costate variables are governed by
\begin{align}
	\dot{p}_x(t) &= -\frac{\partial H}{\partial x}=0 \label{Eq:px}\\
	\dot{p}_y(t) &= -\frac{\partial H}{\partial y}=0 \label{Eq:py}\\
	\dot{p}_\theta(t) &= -\frac{\partial H}{\partial \theta}= p_x(t) \sin\theta(t) - p_y(t) \cos\theta(t)  \label{Eq:ptheta}
\end{align}

It is apparent from Eq.~(\ref{Eq:px}) and Eq.~(\ref{Eq:py}) that $p_x$ and $p_y$ are constant along an optimal trajectory. By integrating Eq.~(\ref{Eq:ptheta}), we have
\begin{align}
	p_\theta(t) =  p_x y(t) - p_y x(t) + c_0
 \label{Eq:ptheta2}
\end{align}
where $c_0$ is a constant. According to the PMP \cite{1987Pontryagin}, we have 
\begin{align}
	\frac{\partial H}{\partial u} =  0
	\label{Eq:H_u}
\end{align}
which can be rewritten explicitly as
\begin{align}
	u(t) =  p_\theta(t), \quad t \in[0,t_f]
	\label{Eq:ut}
\end{align}

Hereafter, we refer to the trajectory that satisfies the necessary conditions Eqs.~(\ref{Eq:px}) - (\ref{Eq:ptheta2}) and Eq.~(\ref{Eq:ut}) as an extremal trajectory.

\section{Real-Time Solution via ANN}
\label{Se:Real-Time}

Let $t_c\in [0,t_f)$ be the current time and let $\boldsymbol{z}_c = (x_c,y_c,\theta_c)$ be the state at $t_c$. Then, the time-to-go is given by
$$t_g = t_f - t_c$$
Furthermore, let us denote by $u^*(t_g,\boldsymbol{z}_c)$ the optimal feedback control at state $\boldsymbol{z}_c$ with a feasible time-to-go $t_g$. By this definition, given an optimal trajectory $\boldsymbol{z}(t)$ of the MECP, if $u(t)$ is the corresponding optimal control, then for any $t\in [0,t_f)$, the following equation holds: 
\begin{equation}
    u^*(t_f-t,\boldsymbol{z}(t)) = u(t)
    \label{Eq:C_u}
\end{equation}
According to the above notations and definitions, addressing the MECP in Problem \ref{Problem1} is equivalent to obtaining the value of the optimal feedback control $u^*(t_g,\boldsymbol{z}_c)$ for any $(t_g,\boldsymbol{z}_c)$ in real time.

\begin{figure}
    \begin{center}
    \includegraphics[width=8.4cm]{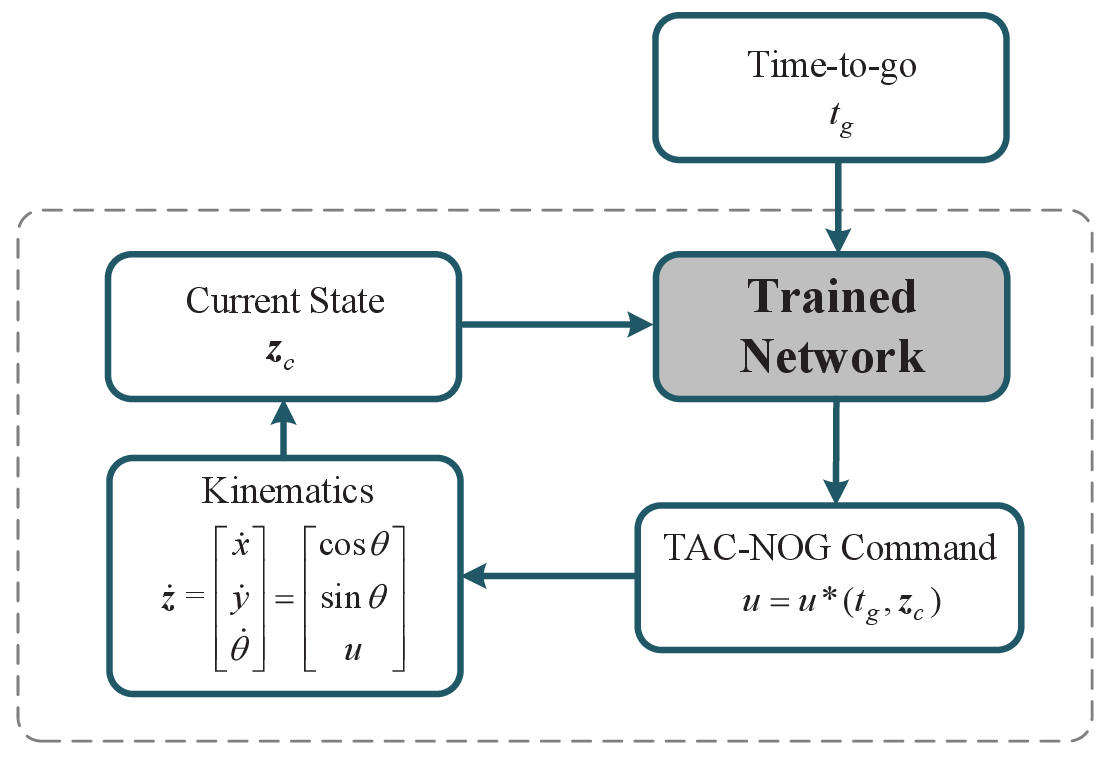}
    \caption{Diagram for the closed-loop guidance via a trained network.}  
    \label{Fig:network}                                 
    \end{center}                                 
\end{figure}

However, it is well known that computing the optimal feedback control $u^*(t_g,\boldsymbol{z}_c)$ in real time is challenging. According to the universal approximation theorem \cite{1989Kurt}, 
if one is able to construct a dataset for the mapping from $(t_g,\boldsymbol{z}_c)$ to $u^*(t_g,\boldsymbol{z}_c)$, a simple ANN can be trained to represent the latter. Note that the output of an ANN is a composition of linear mappings of the input vector. Therefore, for any input vector $(t_g,\boldsymbol{z}_c)$, the optimal feedback control $u^*(t_g,\boldsymbol{z}_c)$ can be generated by the trained ANN within a constant time. If we embed the trained ANN into a closed-loop guidance system, as shown in Fig. ~\ref{Fig:network}, the trained ANN will play the role of generating the TAC-NOG command.

In order to enable an ANN to generate the TAC-NOG command, it is a prerequisite to construct the dataset for the mapping from $(t_g,\boldsymbol{z}_c)$ to $u^*(t_g,\boldsymbol{z}_c)$. In fact, one can use optimization methods, including indirect methods and direct methods, to construct the dataset by solving the MECP in Problem \ref{Problem1} with different initial conditions, as done in \cite{2018Izzo2,2016Sanchez} for soft landing problems. Whereas, both indirect and direct methods suffer from the issue of convergence. In addition, the solution from optimization methods cannot be guaranteed to be at least locally optimal because sufficient conditions for optimality are usually not considered. Thus, if optimization methods are used to construct the dataset, it may include non-optimal solutions. This may further result in the issue that the trained ANN may generate a non-optimal guidance command. In the next section, we shall develop a parameterized approach and establish additional optimality conditions for constructing the dataset of optimal trajectories.

\section{Parameterization of Optimal Trajectories} 
\label{Se:Characterizations}

In this section, we first establish a parameterized differential system so that its solution satisfies the necessary conditions of the MECP. Then, some additional optimality conditions and geometric properties are complemented for the solution of the MECP. By embedding the optimality conditions and geometric properties into the parameterized system, a propagation procedure is presented for generating the dataset for the mapping from $(t_g,\boldsymbol{z}_c)$ to $u^*(t_g,\boldsymbol{z}_c)$. Finally, a scaling invariance property is presented, indicating that it is sufficient to collect the data for the mapping at a fixed time for constructing the dataset.

\subsection{Parameterized System}

Let us define a backward system
\begin{equation}
    \label{Eq:system2}
    \left\{
        \begin{aligned}
		\dot{X} =& -\cos \Theta(t) \\
		\dot{Y} =& -\sin \Theta(t) \\
		\dot{\Theta} =& -[p_x Y(t)-p_y X(t) +c_0]
		\end{aligned}
	\right.
\end{equation}
where $(X,Y)\in \mathbb{R}^2$, $\Theta \in \mathbb{S}$, and the initial condition $[X(0),Y(0),\Theta(0)] $ is fixed as $[0,0,-\pi/2]$. For simplicity, set $\boldsymbol{q}=(p_x,p_y,c_0)$. Then, it is apparent that the solution of Eq.~(\ref{Eq:system2}) is parameterized by $\boldsymbol{q}$ and we shall call it a $\boldsymbol{q}$-parameterized system hereafter. 
Denote $$\boldsymbol{Z}(t,\boldsymbol{q})=(X(t,\boldsymbol{q}),Y(t,\boldsymbol{q}),\Theta(t,\boldsymbol{q}))\in \mathbb{R}^2 \times\mathbb{S}$$
as the solution of the $\boldsymbol{q}$-parameterized system in Eq.~(\ref{Eq:system2}).
It is evident that for any $\boldsymbol{q}$, the solution trajectory $\boldsymbol{Z}(t,\boldsymbol{q})$ of Eq.~(\ref{Eq:system2}) is an extremal trajectory of the MECP in Problem \ref{Problem1}.
Furthermore, set 
$$U(t,\boldsymbol{q}) = p_x Y(t,\boldsymbol{q}) - p_y X(t,\boldsymbol{q})+c_0$$
It is clear that $U(t,\boldsymbol{q})$ is the corresponding extremal control at $\boldsymbol{Z}(t,\boldsymbol{q})$.


\subsection{Optimality Conditions}
\label{Sub:optimality conditions}

In this subsection, we shall present some optimality conditions and geometric properties in addition to the necessary conditions. Let us set
$$\delta(t,\boldsymbol{q}) := \mathrm{det}\left[\frac{\partial \boldsymbol{Z}}{\partial \boldsymbol{q}}(t,\boldsymbol{q}) \right],\ \ \ t\in(0,t_f]$$
as the determinant of the matrix $\frac{\partial \boldsymbol{Z}(t,\boldsymbol{q})}{\partial \boldsymbol{q}}$ for $t\in [0,t_f]$. Then, we have the following result.
\begin{lemma}[Disconjugacy Condition \cite{2016chen}] \label{Le:2}
Given any $\bar{\boldsymbol{q}}\in \mathbb{R}^3$ and any $\bar{t}>0$, if $\delta(t,\bar{\boldsymbol{q}})\neq 0$ for $t\in (0,\bar{t})$, then we have that $\boldsymbol{Z}(t,\bar{\boldsymbol{q}})$ for $t\in [0,\bar{t}]$ is at least a locally optimal trajectory. On the other hand, if $\delta(t,\bar{\boldsymbol{q}}) = 0$ for a time in $(0,\bar{t}]$, then $\boldsymbol{Z}(t,\bar{\boldsymbol{q}})$ for $t\in [0,\bar{t}]$ loses its local optimum.
\end{lemma}

Lemma \ref{Le:2} is a direct result of Theorem 2 in \cite{2016chen}, and it presents a sufficient condition for optimality.
By the following lemma, we shall present another optimality condition by analyzing the geometric property of the extremal trajectory.

\begin{lemma}\label{Le:3}
Given any $\bar{\boldsymbol{q}}\in \mathbb{R}^3$, if there are two different instants $t_1\in(0,t_f)$ and $t_2\in (0,t_f)$ so that the velocity vectors $[\cos\Theta(t_1,\bar{\boldsymbol{q}}), \sin\Theta(t_1,\bar{\boldsymbol{q}})]$ and $[\cos\Theta(t_2,\bar{\boldsymbol{q}}),\sin\Theta(t_2,\bar{\boldsymbol{q}})]$ are colinear, i.e.,
\begin{align}
\label{Eq:co-linear}
    \frac{Y(t_1,\bar{\boldsymbol{q}})-Y(t_2,\bar{\boldsymbol{q}})}{X(t_1,\bar{\boldsymbol{q}})-X(t_2,\bar{\boldsymbol{q}})} = \tan\Theta(t_1,\bar{\boldsymbol{q}}) = \tan\Theta(t_2,\bar{\boldsymbol{q}})
\end{align}
then the extremal trajectory $\boldsymbol{Z}(t,\bar{\boldsymbol{q}})$ for $t\in[0,t_f]$ is not optimal.
\end{lemma}

The proof of Lemma \ref{Le:3} is postponed to Appendix \ref{App:C}.
Lemma \ref{Le:2} and Lemma \ref{Le:3} introduce two new optimality conditions for extremal trajectories.  
Therefore, the trajectory
obtained by the $\boldsymbol{q}$-parameterized system that satisfies all the necessity conditions (Eqs.~(\ref{Eq:px}) - (\ref{Eq:ptheta2}) and Eq.~(\ref{Eq:ut})) along with the supplementary conditions (Lemma \ref{Le:2} - \ref{Le:3}) is at least a locally optimal solution.

\subsection{Numerical Procedure for Generating the Dataset}

According to the developments in the previous subsections, for any given $\boldsymbol{q}$, a locally optimal trajectory can be obtained by integrating the $\boldsymbol{q}$-parameterized system in Eq.~(\ref{Eq:system2}).
Therefore, by randomly choosing values for $\boldsymbol{q}$, we can obtain a large number of optimal trajectories. One is able to collect, from these optimal trajectories, the dataset for the mapping from $(t_g,\boldsymbol{z}_c)$ to $u^*(t_g,\boldsymbol{z}_c)$. The following lemma shall show that collecting the data for the mapping at a fixed instant instead of the whole trajectories is sufficient for training an ANN to generate the TAC-NOG command.
\begin{lemma}
\label{Le:zoom}
Given a state $(x_c,y_c,\theta_c)\in \mathbb{R}^2 \times\mathbb{S}$ and a time-to-go $t_g>0$, assume that there exists an optimal trajectory $(x(t),y(t),\theta(t))$ for $t\in [0,t_g]$ so that  $(x(0),y(0),\theta(0)) = (x_c,y_c,\theta_c)$ and $(x(t_g),y(t_g),\theta(t_g)) = (0,0,-\pi/2)$. Then, for any $T\in \mathbb{R}_+$ we have
\begin{align}
\label{Eq:transform2}
u^*(t_g,x_c,y_c,\theta_c) = \frac{T}{t_g} u^*(T,\frac{T}{t_g}x_c,\frac{T}{t_g}y_c,\theta_c)
\end{align}
\end{lemma}
The proof of this lemma is postponed to Appendix \ref{App:C}.
Let us denote by $\mathcal{Z}_T\in \mathbb{R}^2\times \mathbb{S}$ the set of state $(x,y,\theta)$ that can be controlled to the final condition $(0,0,-\pi/2)$ with the duration of $T$. In addition, let us denote by $N(\boldsymbol{z})$ the mapping from $(T,\boldsymbol{z})$ with $\boldsymbol{z}\in \mathcal{Z}_T$ to $u^*(T,\boldsymbol{z})$. 
Then, Lemma \ref{Le:zoom} indicates that for any $(t_g,x_c,y_c,\theta_c)$ with $t_g \neq T$, we can find a state $(x,y,\theta)\in \mathcal{Z}_T$ with $x=Tx_c/t_g$, $y=Ty_c/t_g$, and $\theta= \theta_c$ so that
\begin{align}
\label{Eq:chen_jia}
u^*(t_g,\boldsymbol{z}_c) = \frac{T}{t_g} N(\boldsymbol{z})  
\end{align}
where $\boldsymbol{z}_c = (x_c,y_c,\theta_c)$ and $\boldsymbol{z}=(x,y,z) \in \mathcal{Z}_T$.
As a consequence, we just need to use a neural network to approximate the mapping $N(\boldsymbol{z}): \mathcal{Z}_T\rightarrow \mathbb{R}$.
\begin{lemma}
\label{Le:chen}
For any state $\boldsymbol{z}\in \mathcal{Z}_T$, there exists a $\boldsymbol{q}\in \mathbb{R}^3$ so that $\boldsymbol{z} = \boldsymbol{Z}(T,\boldsymbol{q})$  and $N(\boldsymbol{z}) = U(T,\boldsymbol{q})$. On the other hand, for any $\boldsymbol{q}\in \mathbb{R}^3$ we have $\boldsymbol{Z}(T,\boldsymbol{q})\in \mathcal{Z}_T$.
\end{lemma} 
The proof of this lemma is postponed to Appendix \ref{App:C}. Thanks to Lemma \ref{Le:chen}, we are able to construct the dataset $\mathcal{D}$ for the mapping $N(\boldsymbol{z}): \mathcal{Z}_T\rightarrow \mathbb{R}$ by propagating the parameterized system, as shown in Algorithm \ref{algorithm_1}.
\begin{algorithm} 
	\caption{Generation of the Dataset}   
	 \label{algorithm_1}       
    \begin{algorithmic}[1] 
    \REQUIRE positive numbers $p_{max}$ and $T$
    \STATE $\mathcal{D}\gets \varnothing$
    \FOR {$p_x = -p_{max}$ to $p_{max}$}
        \FOR{$p_y = -p_{max}$ to $p_{max}$}
            \FOR{$c_0 = -p_{max}$ to $p_{max}$}
                \STATE Let $\boldsymbol{q}=(p_x,p_y,c_0)$ and propagate the parameterized system in Eq.~(\ref{Eq:system2}) with the initial condition $(0,0,-\pi/2)$ to get the trajectory $Z(t,\boldsymbol{q})$ and the control $U(t,\boldsymbol{q})$  for $t\in [0,T]$
                \IF {$\delta(t,\boldsymbol{q}) \neq 0$ for $t\in (0,T]$ and if Eq.~(\ref{Eq:co-linear}) is not met for $t\in (0,T]$}
                    \STATE $\mathcal{D} \leftarrow \mathcal{D} \cup \{\boldsymbol{Z}(T,\boldsymbol{q}),U(T,\boldsymbol{q})\}$\;
               \ENDIF
           \ENDFOR
        \ENDFOR
   \ENDFOR
   \end{algorithmic} 
\end{algorithm}

A simple ANN can be trained by the dataset $\mathcal{D}$ generated in Algorithm \ref{algorithm_1} to represent the mapping $N(\boldsymbol{z}): \mathcal{Z}_T\rightarrow \mathbb{R}$. Then, one can combine the trained ANN with the transformation in Eq.~(\ref{Eq:transform2}) to generate the TAC-NOG command within a constant time. In the next section, we shall present some numerical examples to illustrate how the trained ANN is used to generate the TAC-NOG command. 

\section{Numerical Simulations}
\label{Se:Simulations}

Set the parameters $p_{max}$ and $T$ in Algorithm \ref{algorithm_1} as 30 and 1.5, respectively. Meanwhile, the step size for the iterations in Algorithm \ref{algorithm_1} is set as 0.05. Then, the number of elements in dataset $\mathcal{D}$ is no more than $1.846\times10^8 $. An ANN with two hidden layers, each of which contains 20 neurons, is trained to represent the mapping $N(\boldsymbol{z}): \mathcal{Z}_T\rightarrow \mathbb{R}$. For notational simplicity, we denote by $\bar{N}(\boldsymbol{z})$ as the trained ANN. Then, according to Eq.~(\ref{Eq:chen_jia}), for any $(t_g,\boldsymbol{z}_c)$ with $\boldsymbol{z}_c=(x_c,y_c,\theta_c)$ we have
$$u^*(t_g,\boldsymbol{z}_c) \approx \frac{T}{t_g}\bar{N}(\frac{T}{t_g}x_c,\frac{T}{t_g}y_c,\theta_c)$$

Given any $\boldsymbol{z}\in \mathcal{Z}_T$ as input, the output of the trained network, $\bar{N}(\boldsymbol{z})$, can be obtained in around 6.72 $\rm{\mu s}$ using a laptop with Intel Core i5-10400 CPU @2.90 GHz. Therefore, this trained network can be used to approximate the optimal control command in real time. 
To further demonstrate the merits of the trained network in generating TAC-NOG command, we will present some examples below by comparing with optimization methods and existing guidance methods.

\subsection{Comparisons with Optimization Methods}\label{SubSe:global optimality}

In this subsection, we shall compare the performance of the TAC-NOG law developed in this paper with a shooting method.

\begin{figure}
    \begin{center}
    \includegraphics[width=0.4\textwidth]{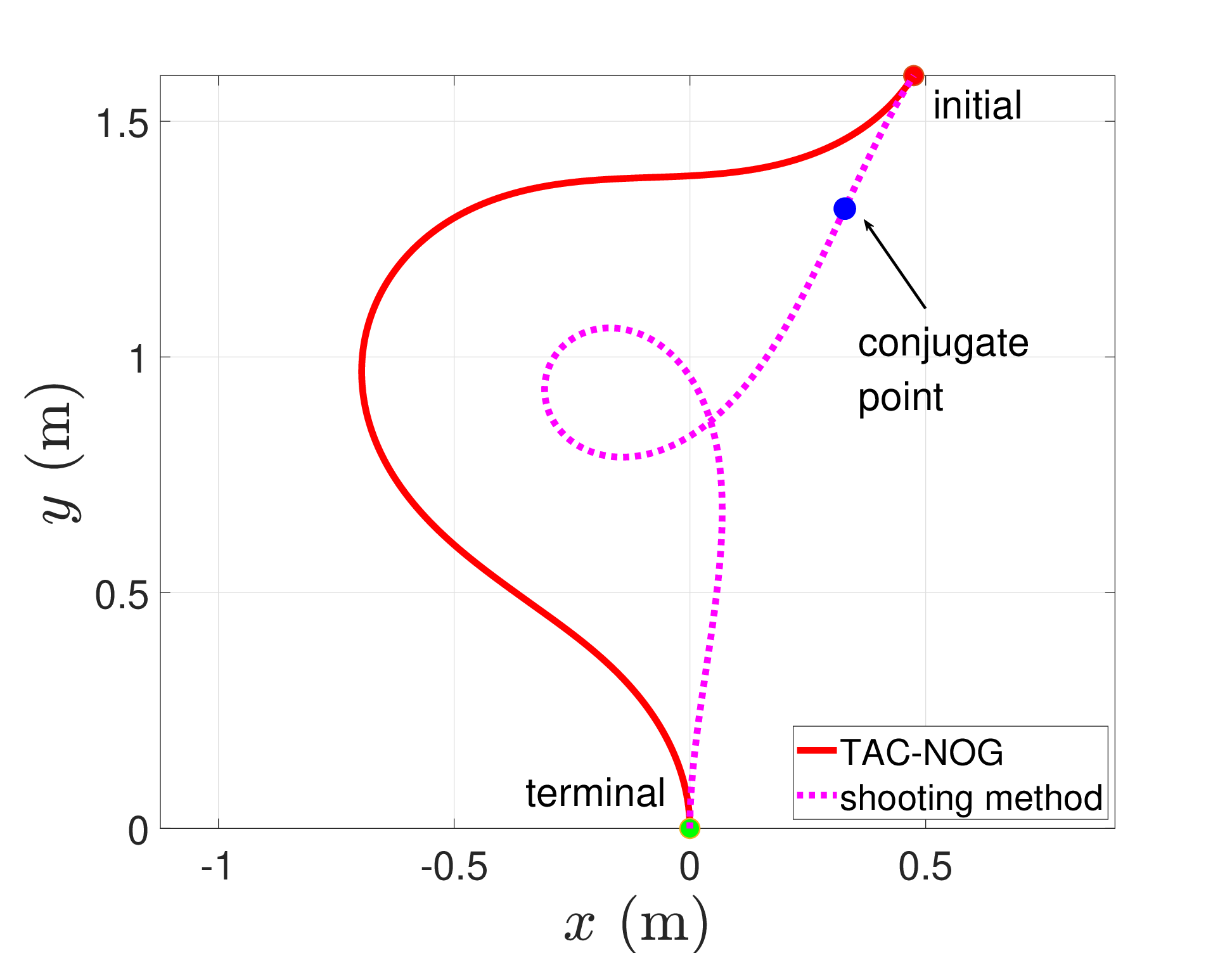}    
    \caption{Trajectories generated by TAC-NOG law and shooting method.}  
    \label{Fig:local_trajectory}                               
    \end{center}                                 
\end{figure}

\begin{figure}
    \centering
    \subfigure[Controls]{
    \label{Fig:local_control}
        \centering
        \includegraphics[width=0.45\textwidth]{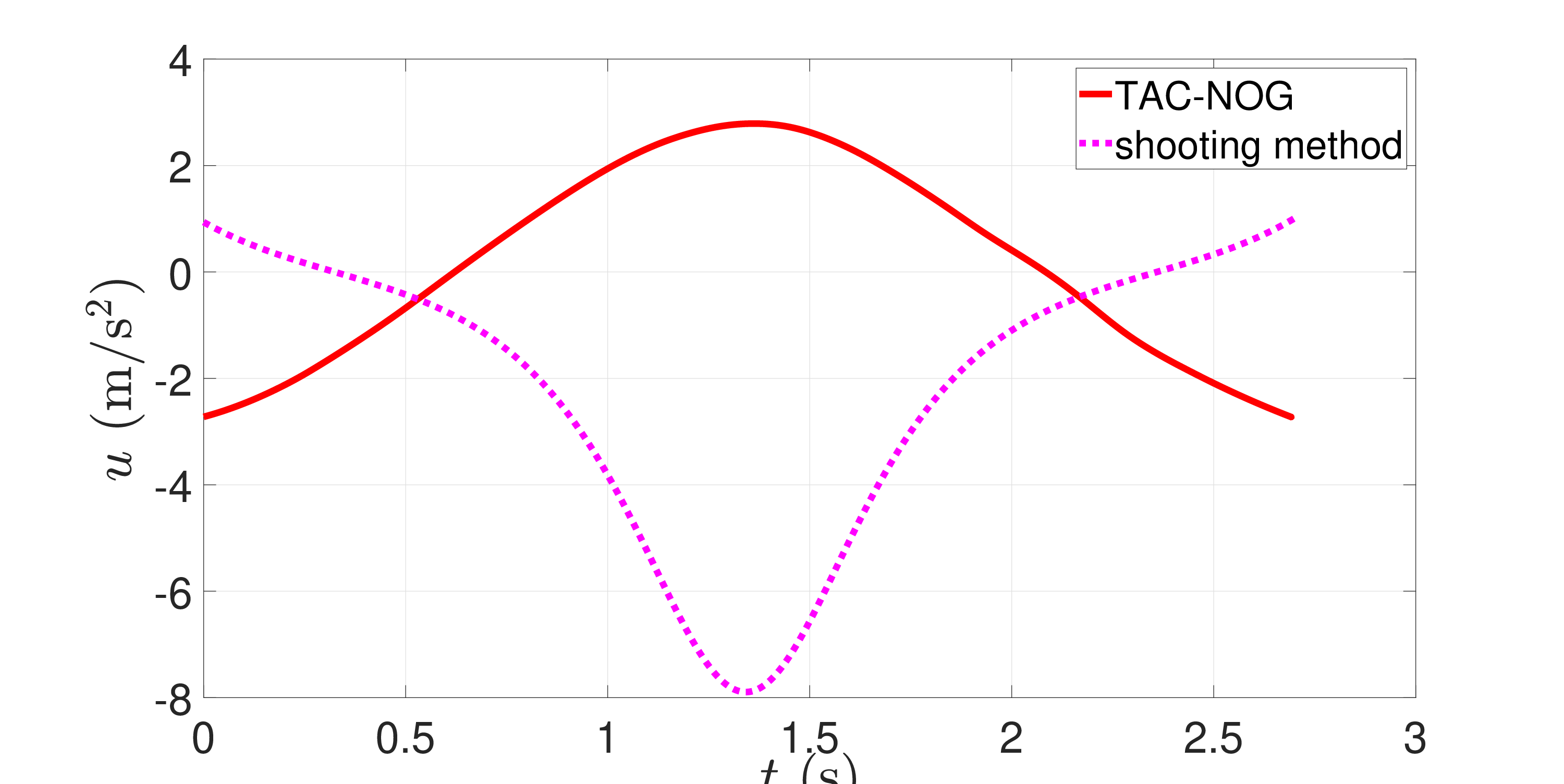} }
        \hfill
        \subfigure[Heading angles]{
    \label{Fig:local_angle}
        \centering
        \includegraphics[width=0.45\textwidth]{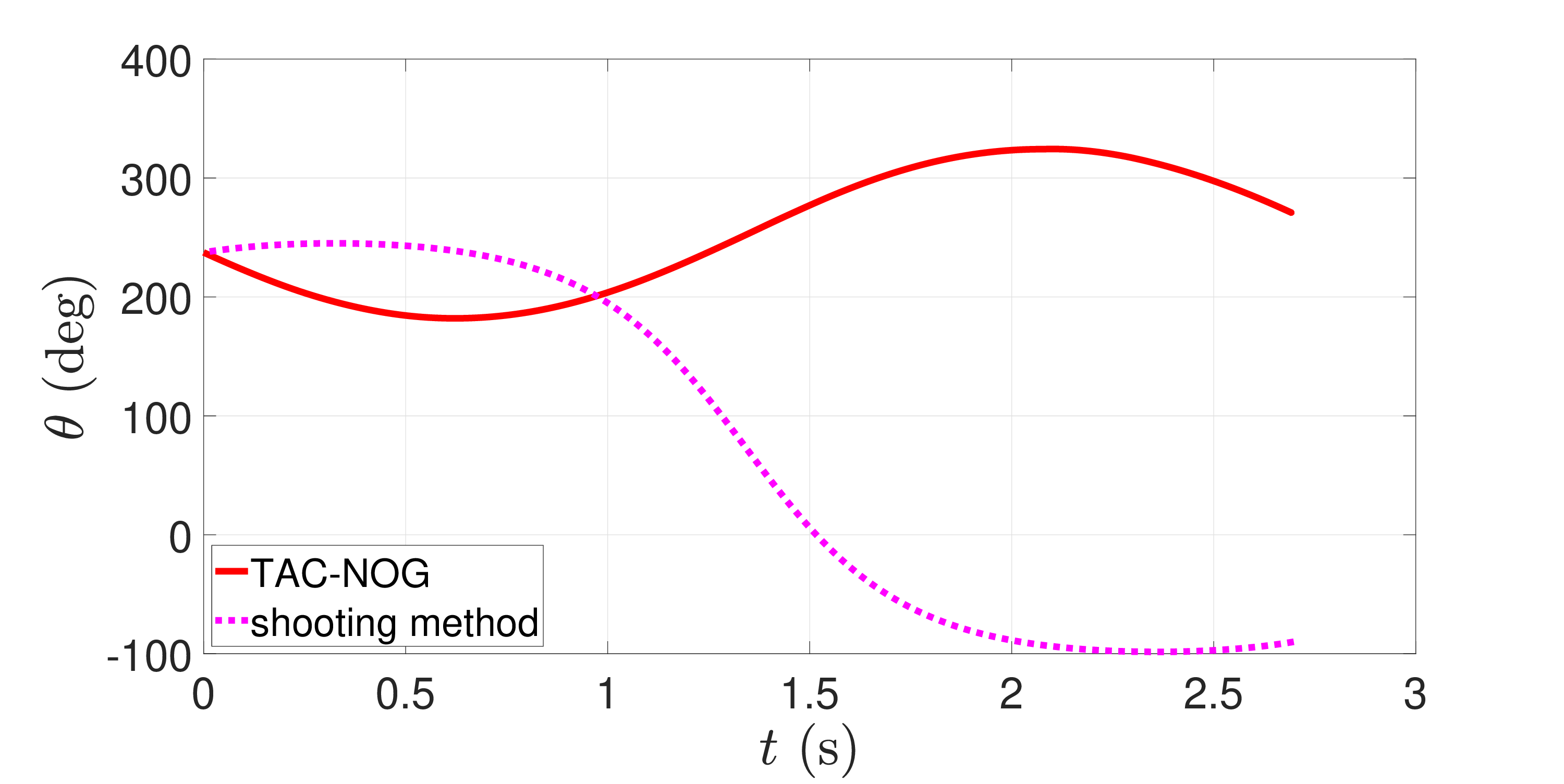}
    }
    
    \caption{Control profiles and heading angles related to TAC-NOG law and shooting method.}
    \label{Fig:local}
\end{figure}

Consider the initial position of the pursuer to be (0.4748, 1.5968) m with an initial heading angle of 237.4 deg. The speed of the pursuer is 1 m/s. 
The target is located at the origin (0, 0) m.
The pursuer is expected to intercept the target with an impact angle of -90 deg and a desired impact time of 2.7 s.
Trajectories generated by the TAC-NOG law and the shooting method are presented in Fig.~\ref{Fig:local_trajectory}. The profiles of control commands and look angles along these two trajectories are demonstrated in Fig.~\ref{Fig:local}.

\begin{figure}
    \begin{center}
    \includegraphics[width=0.45\textwidth]{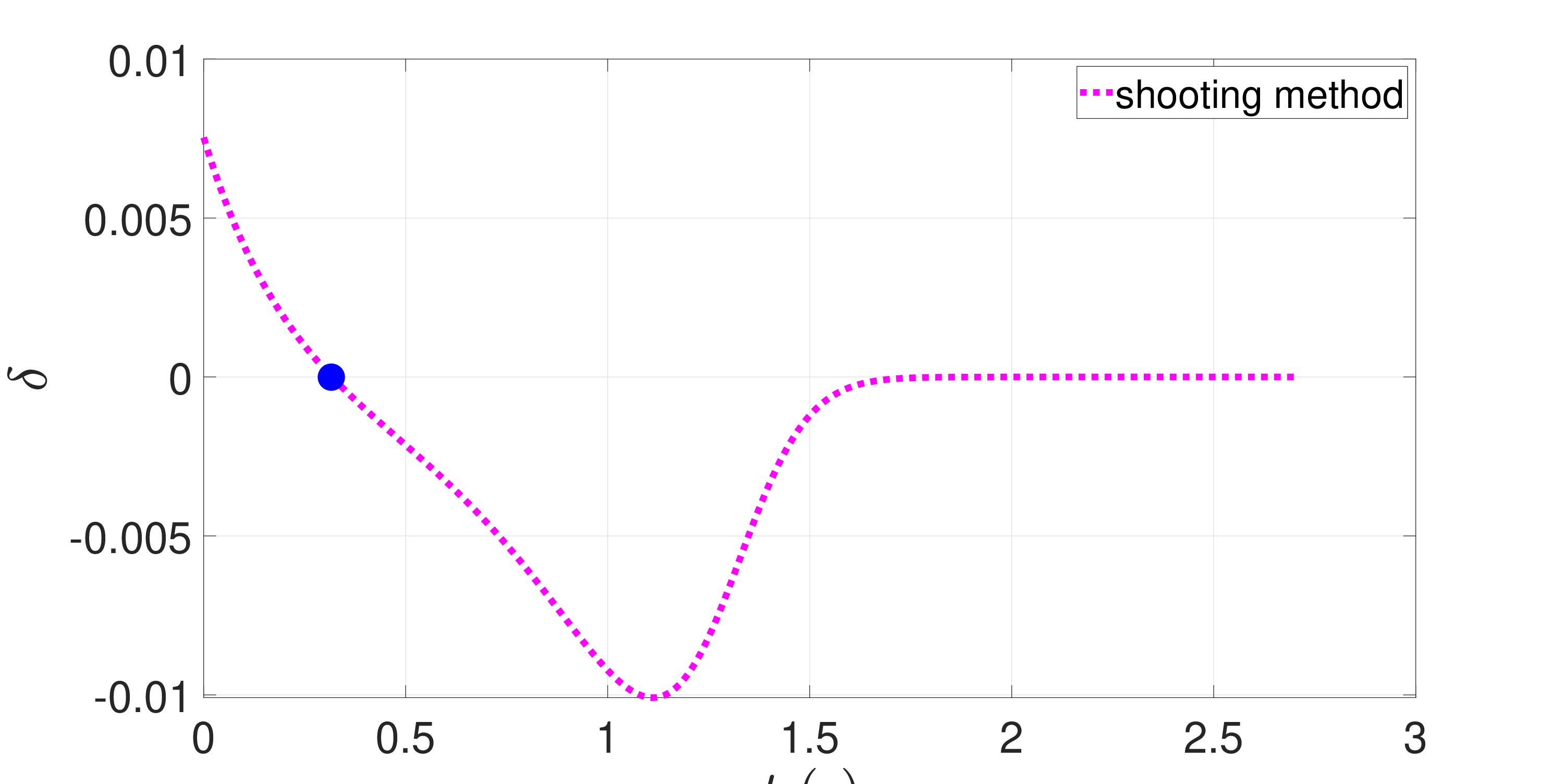}    
    \caption{The time history of $\delta(t,\boldsymbol{q})$ along the trajectory of shooting method.}  
    \label{Fig:conjugate time}                               
    \end{center}                                 
\end{figure}

From Fig.~\ref{Fig:local_trajectory} and Fig.~\ref{Fig:local}, it is evident that both methods allow the pursuer to reach the target with the expected time and terminal heading angle. However, the trajectories diverge significantly; see Fig.~\ref{Fig:local_trajectory}. Note that the shooting method only considers necessary conditions from PMP, which cannot guarantee the solution to be at least locally optimal. 
We can see from Fig.~\ref{Fig:conjugate time} that the function $\delta(t,\boldsymbol{q})$ changes sign at a point along the trajectory of the shooting method. This indicates that the trajectory of the shooting method is not optimal; see Lemma \ref{Le:2}. 
In fact, the control effort required by the TAC-NOG law is 4.5567 $\rm{m^2/s^3}$, while that related to the shooting method is up to 16.1776 $\rm{m^2/s^3}$. This coincides with the control profiles in Fig.~\ref{Fig:local_control}. It is worth mentioning that all the optimality conditions in Lemma $\ref{Le:2}$ and Lemma $\ref{Le:3}$ are satisfied along with the solution generated by TAC-NOG law. Thus, the solution generated by the TAC-NOG law is at least locally optimal.

\subsection{Comparisons with Existing Guidance Laws}

In this subsection, we will compare the performance of the developed TAC-NOG law with two existing guidance laws: the BPNG law \cite{2013Zhang} and the SMC-based guidance law \cite{2019Chenxiaot}. The pursuer’s normal acceleration is limited within $5g$, where  $g=9.8~\rm{m/s^2}$ is the gravitational acceleration.

\textbf{Case A: Comparisons with BPNG Law}

Consider the initial position of the pursuer situated at (-10000, 1000) m, with an initial heading angle of 60 deg. The target is located at (500, 0) m. 
Set the speed of the pursuer as 250 m/s.
The pursuer aims to intercept the target with an impact angle of 10 deg, and we consider four different impact times: 45 s, 50 s, 55 s, and 60 s, for simulation.

\begin{figure}
    \begin{center}
    \includegraphics[width=0.4\textwidth]{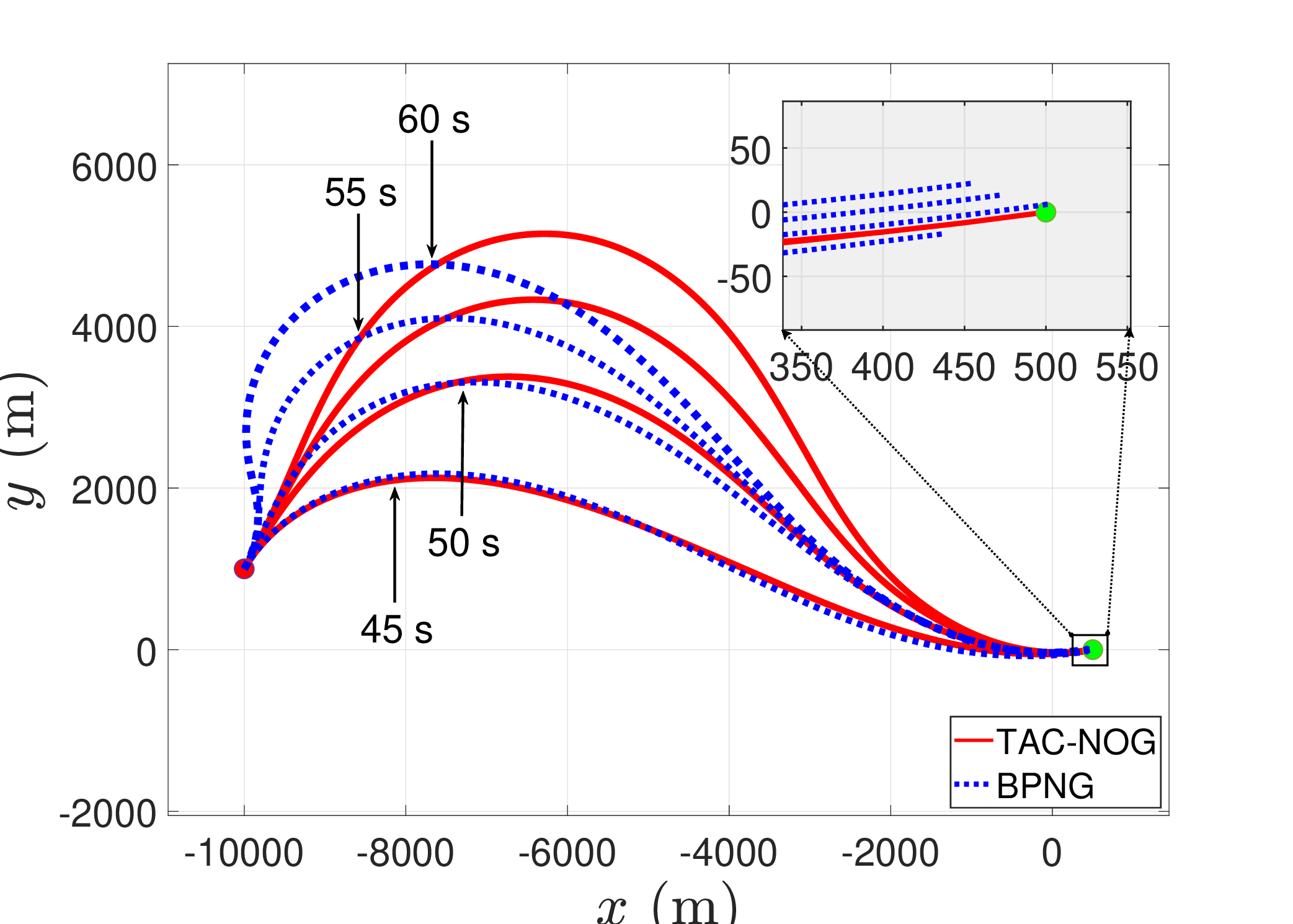}    
    \caption{Case A: trajectories generated by TAC-NOG law and BPN law with different impact times.}  
    \label{Fig:different_time_1}                               
    \end{center}                                 
\end{figure}

\begin{figure}
    \centering
    \subfigure[$t_f=45$ s]{
        \centering
        \includegraphics[width=0.22\textwidth]{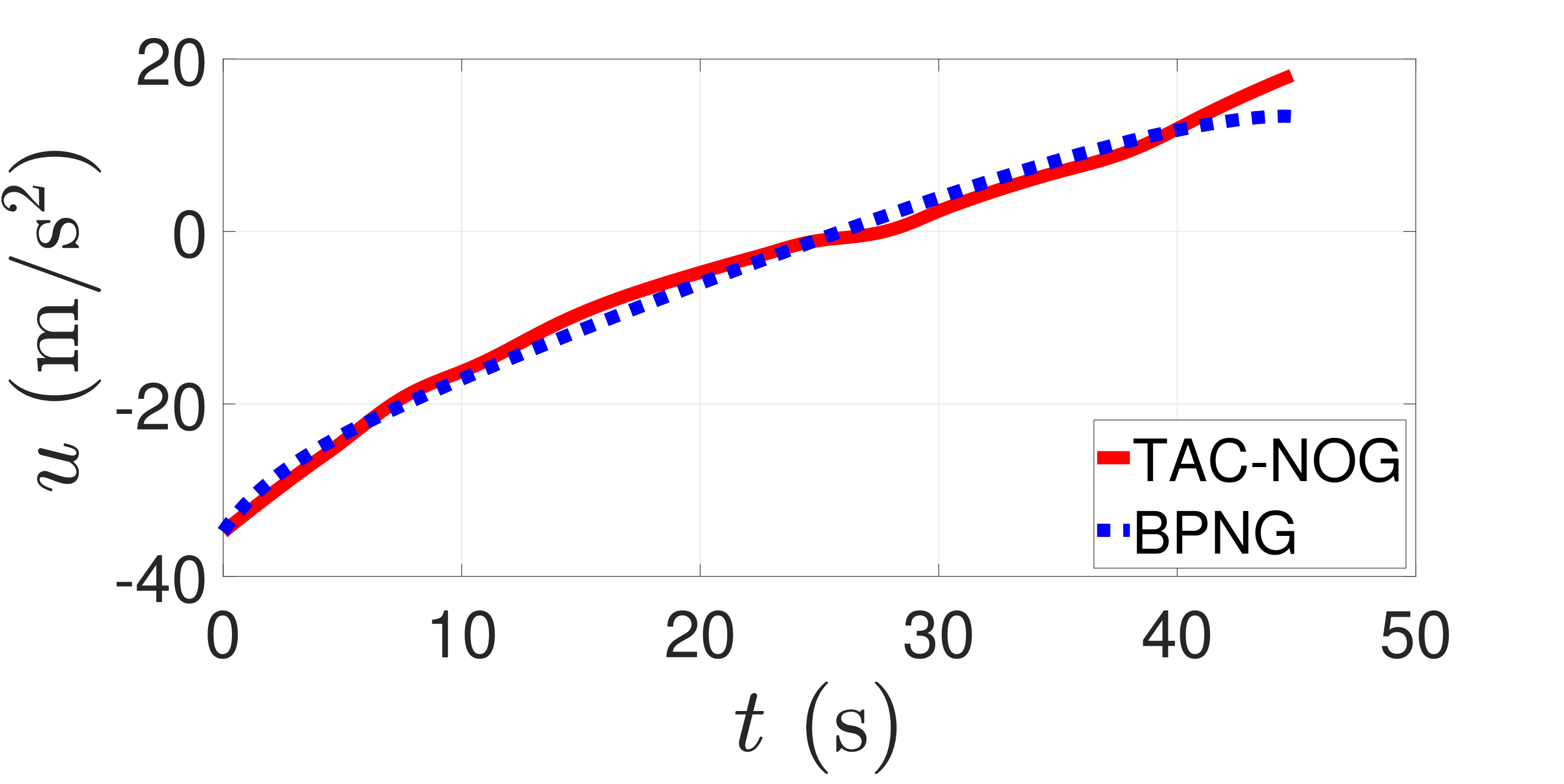} }
        \hspace{-3mm}
    \subfigure[$t_f=50$ s]{
        \centering
        \includegraphics[width=0.22\textwidth]{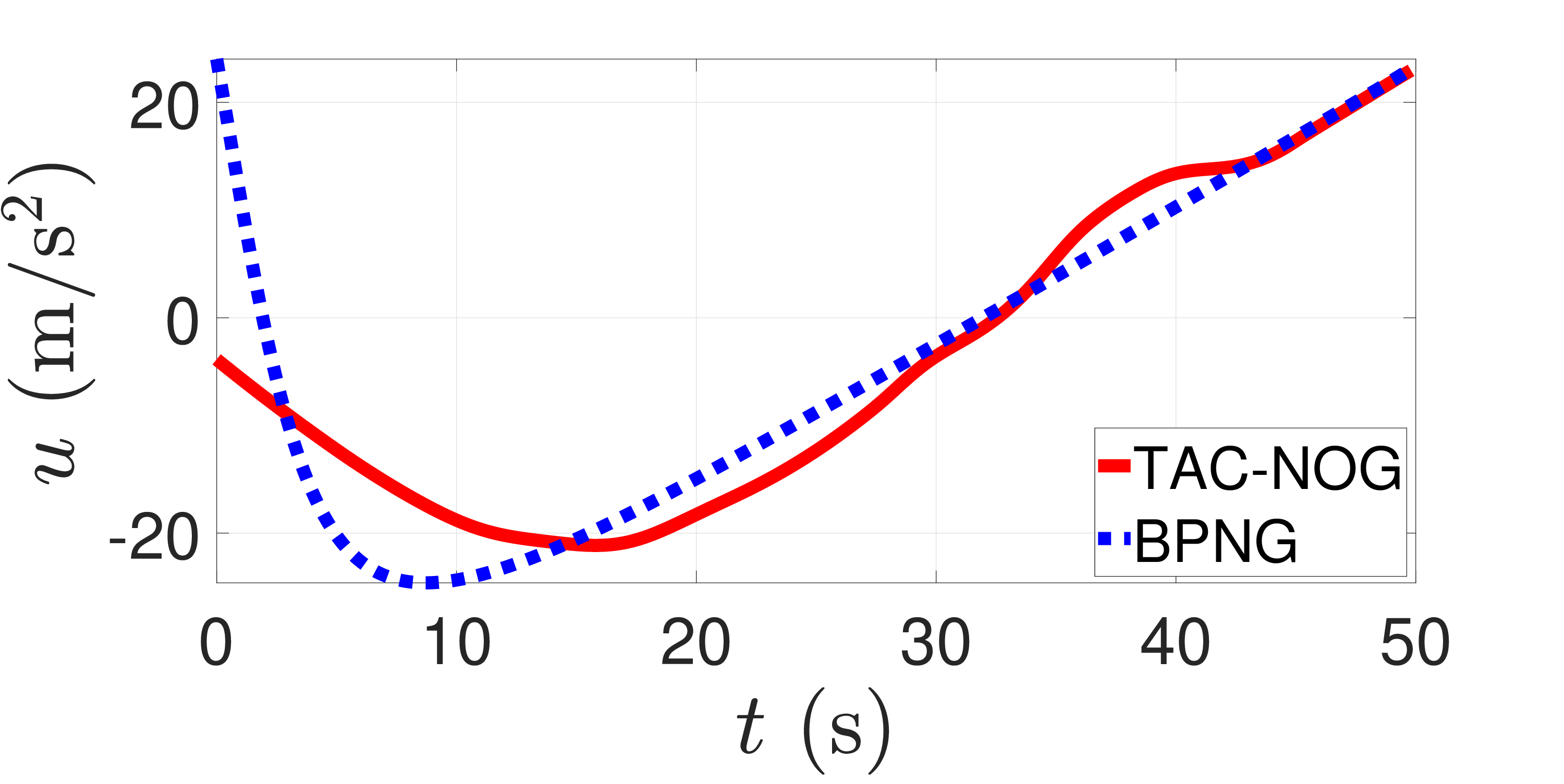} }
        \hfill
    \subfigure[$t_f=55$ s]{
        \centering
        \includegraphics[width=0.22\textwidth]{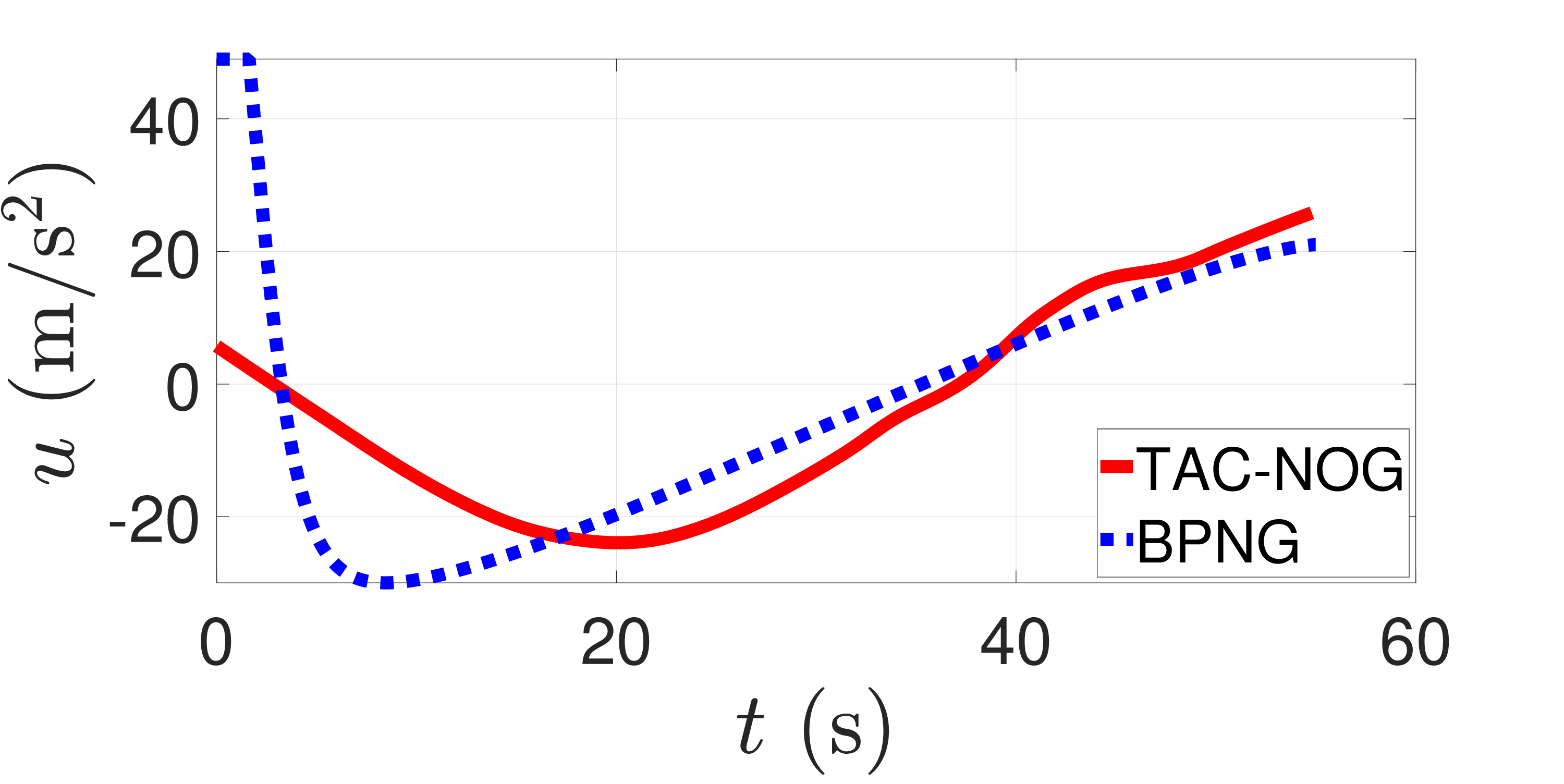} }
        \hspace{-3mm}
    \subfigure[$t_f=60$ s]{
        \centering
        \includegraphics[width=0.22\textwidth]{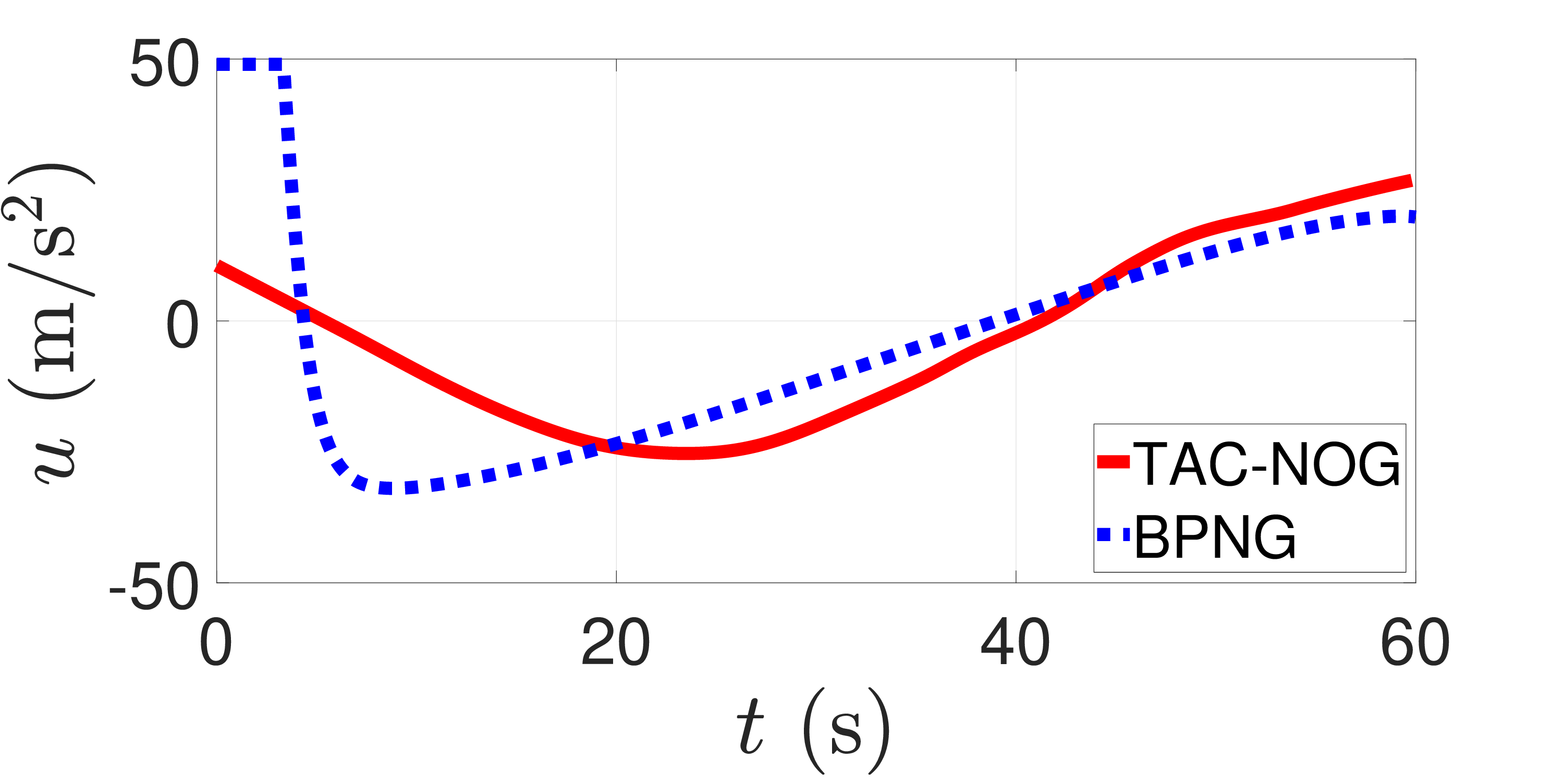}            
    }
    \caption{Case A: control profiles related to TAC-NOG law and BPNG law.}
    \label{Fig:different_time_2}
\end{figure}

\begin{figure}
    \begin{center}
    \includegraphics[width=0.45\textwidth]{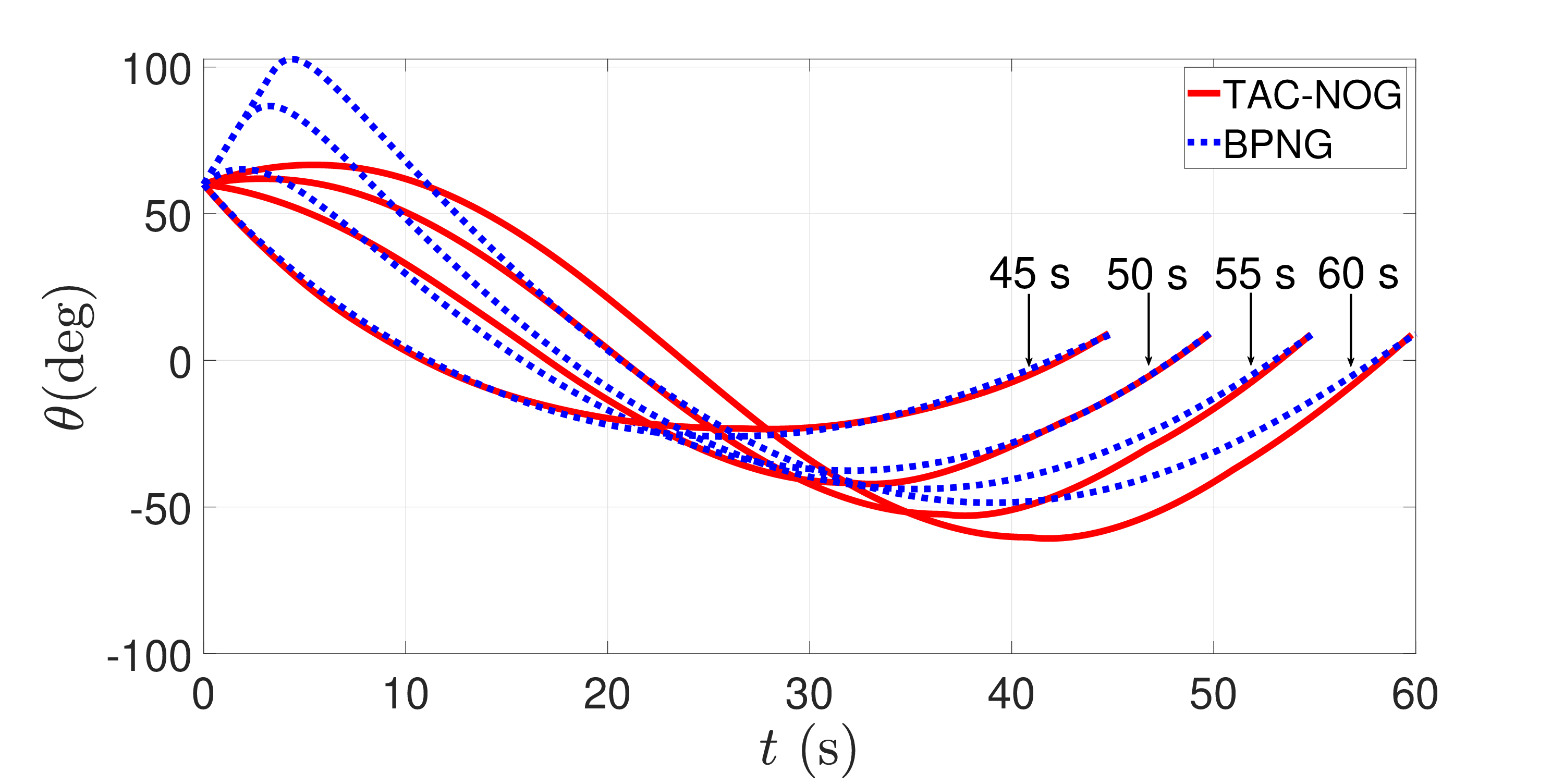}    
    \caption{Case A: heading angles related to TAC-NOG law and BPNG law.}  
    \label{Fig:different_time_3}                            
    \end{center}                                 
\end{figure}

The trajectories generated by the TAC-NOG law and the BPNG law are shown by solid curves and dotted curves in Fig.~\ref{Fig:different_time_1}, respectively. The time histories of corresponding controls and heading angles are presented in Fig.~\ref{Fig:different_time_2} and Fig.~\ref{Fig:different_time_3}, respectively. Though the expected impact angle and impact time can be precisely achieved by the TAC-NOG law and the BPNG law, notable differences can be seen in Fig.~\ref{Fig:different_time_1} and Fig.~\ref{Fig:different_time_2} in terms of trajectories and control commands.

From Fig.~\ref{Fig:different_time_1}, it is evident that the errors of the terminal position associated with the BPNG law are larger than those associated with the TAC-NOG law. In fact, the largest error among the four trajectories generated by the BPNG is up to 64.7 m, whereas it can be reduced to 2.9 m by the TAC-NOG law. 
Additionally, we can see from Fig.~\ref{Fig:different_time_2} that the controls required by the BPNG are relatively larger during the early stage of the pursuing process. This indicates that the TAC-NOG law needs less control effort, as shown by the data in Table \ref{Tab:table1}, where $J_{TAC-NOG}$ and $J_{BPNG}$ denote the control efforts required by TAC-NOG law and BPNG law, respectively. This numerical result coincides with the theoretical development as the TAC-NOG law is an optimal guidance law.

\begin{table}[!hbt]
\centering
\caption{Case A: the control effort required by TAC-NOG law and BPNG law for different impact times.}
\label{Tab:table1}
\begin{tabular}{c|cc}
\hline
$t_f~(\rm{s})$ & $J_{TAC-NOG}~(\rm{m^2/s^3})$ & $J_{BPNG}~(\rm{m^2/s^3})$ \\\hline
45 & $4.578\times 10^3$ & $4.581\times 10^3$ \\
50 & $5.435\times 10^3$ & $6.009\times 10^3$ \\
55 & $7.174\times 10^3$ & $1.077\times 10^4$ \\
60 & $8.629\times 10^3$ & $1.427\times 10^4$\\
\hline
\end{tabular}
\end{table}

\textbf{Case B: Comparisons with BPNG Law and SMC-based Law}

Consider a scenario involving four pursuers intercepting a target at the same impact time from different impact angles. The target is located at (1000, 5000) m.
The initial state $(x_0,y_0,\theta_0)$ for each pursuer is given as

\begin{align}
    \nonumber
    Pursuer \ \#1:~& \boldsymbol{z}(0)=(-15000~\rm{m}, -10000~\rm{m},-110~\rm{deg})\\
\nonumber
    Pursuer \ \#2:~& \boldsymbol{z}(0)=(5000~\rm{m},-20000~\rm{m},45~\rm{deg})\\
    \nonumber
    Pursuer \ \#3: ~&\boldsymbol{z}(0)=(18000~\rm{m}, 20000~\rm{m}, 130~\rm{deg})\\
    \nonumber
    Pursuer \ \#4: ~&\boldsymbol{z}(0)=(-10000~\rm{m}, 10000~\rm{m},-130~\rm{deg})
\end{align}

For notational simplicity, we denote pursuer $\#i$ by $P_i$, denote the speed of $P_i$ by $V_i$, and denote the impact angle of $P_i$ by $\theta_{if}$. Set the values of $V_1$, $V_2$, $V_3$, and $V_4$ as 250 m/s, 250 m/s, 300 m/s, and 200 m/s, respectively, and set the values of $\theta_{1f}$, $\theta_{2f}$, $\theta_{3f}$, and $\theta_{4f}$ as 90 deg, 180 deg, 270 deg, and 0 deg, respectively. The desired impact times for all four pursuers are set the same as 120 s. 

\begin{figure}
    \begin{center}
    \includegraphics[width=0.4\textwidth]{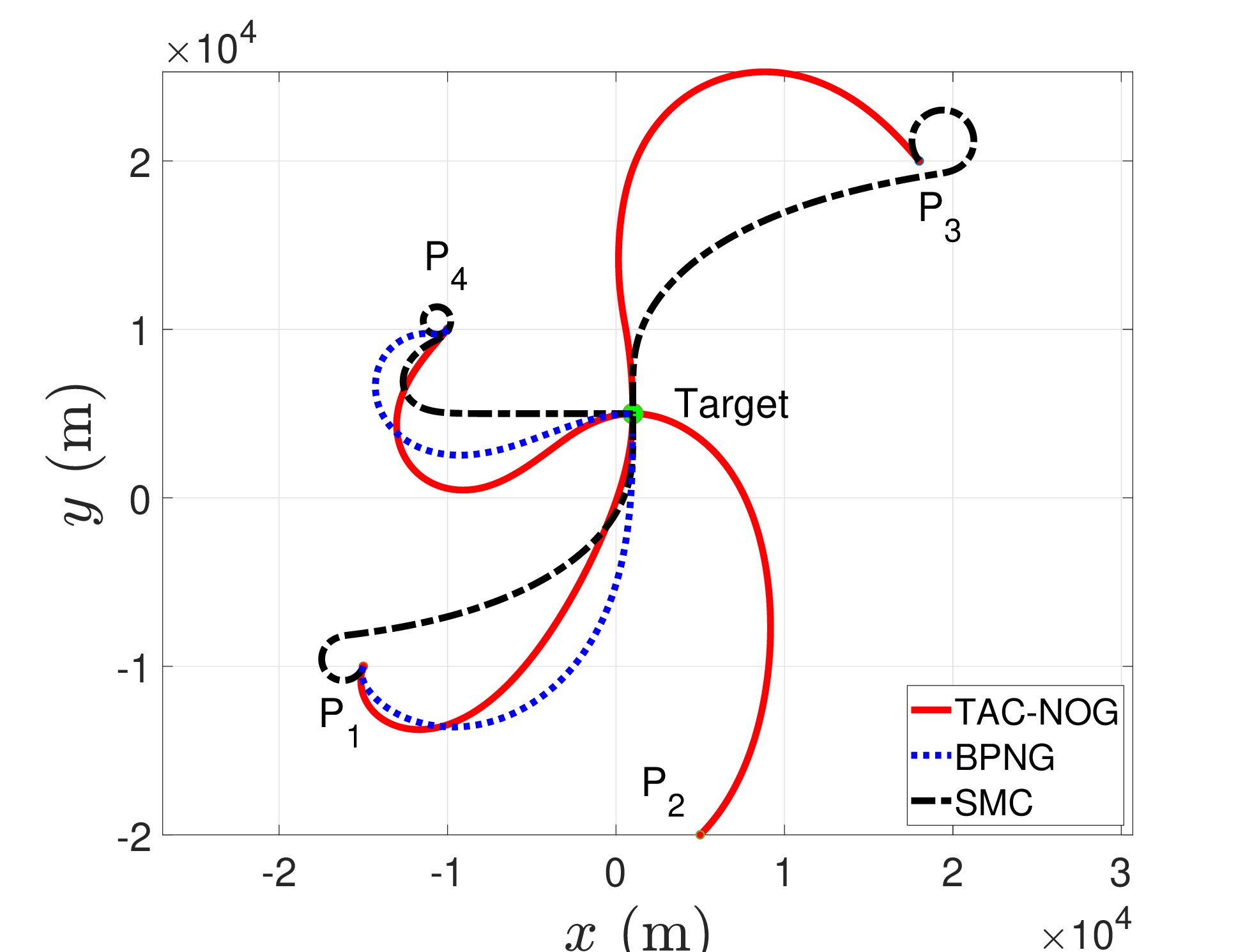}    
    \caption{Case B: trajectories generated by three guidance laws with different impact angles.}  
    \label{Fig:different angle}                                 
    \end{center}                                 
\end{figure}

The trajectories generated by the TAC-NOG law, the BPNG law and the SMC-based law \cite{2019Chenxiaot} are presented by the solid curves, dotted curves and dashed curves in Fig.~\ref{Fig:different angle}, respectively. The guidance command for each pursuer is computed individually. We can observe from Fig.~\ref{Fig:different angle} that $P_1$ and $P_4$ can be controlled to the target with the expected impact angle and impact time using any of the three guidance laws. However, the BPNG law operates under the assumption of small angles, limiting its applicability to a narrow range of impact angles. Consequently, for $P_2$ and $P_3$, the BPNG law is not applicable.
When using the SMC-based guidance law, $P_2$ is unable to reach the target with the expected impact angle and impact time. This occurs because the pursuer is initially guided to an inappropriate sliding surface at 12.772 s. Subsequently, the guidance law switches to the PN guidance law, for the remainder of the interception, resulting in unsatisfactory impact time and impact angle.

\begin{figure}
    \centering
    \subfigure[Pursuer $\#1$]{
        \centering
        \includegraphics[width=0.22\textwidth]{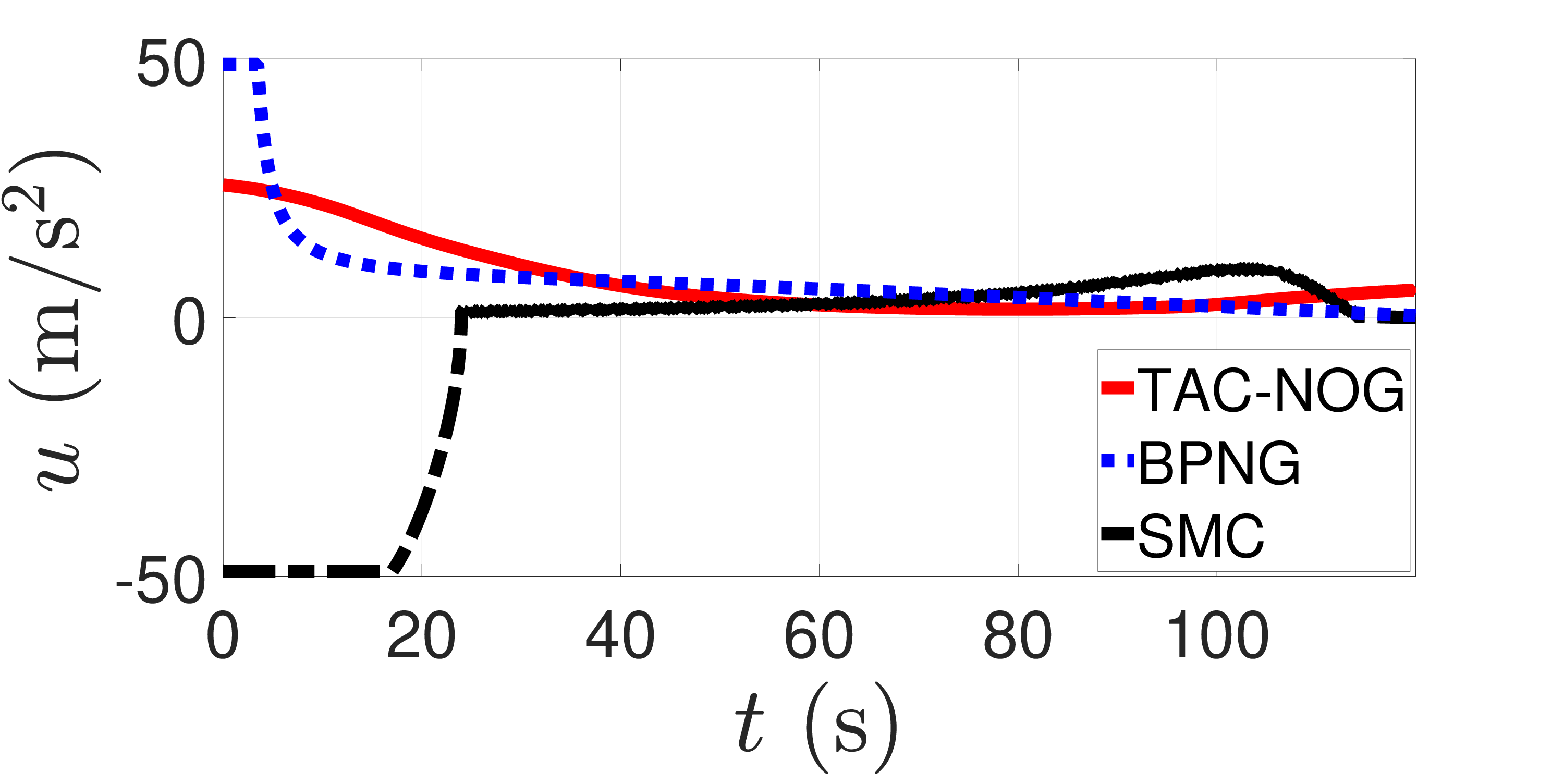} }
        \hspace{-3mm}
    \subfigure[Pursuer $\#2$]{
        \centering
        \includegraphics[width=0.22\textwidth]{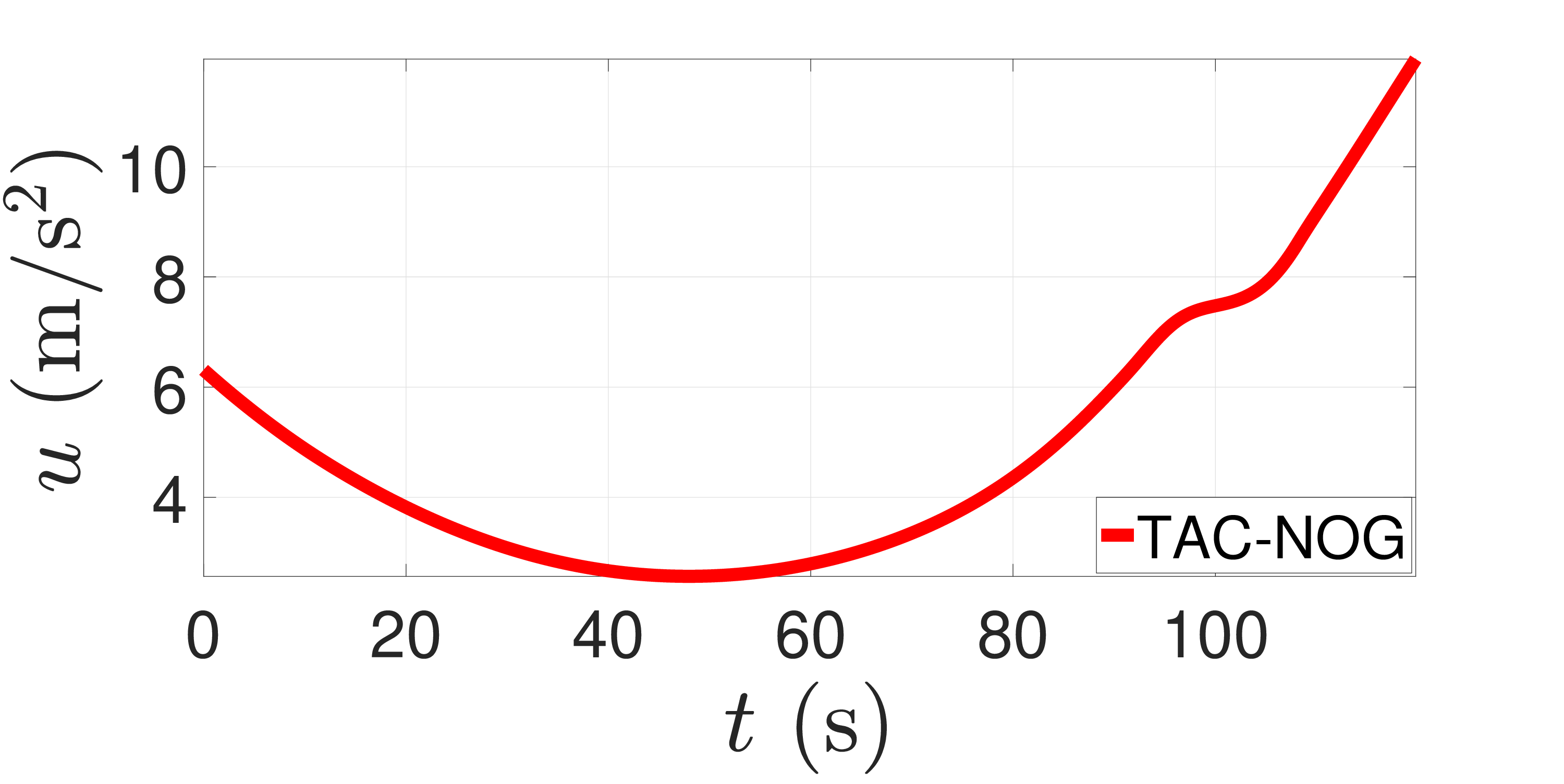} }
        \hfill
    \subfigure[Pursuer $\#3$]{
        \centering
        \includegraphics[width=0.22\textwidth]{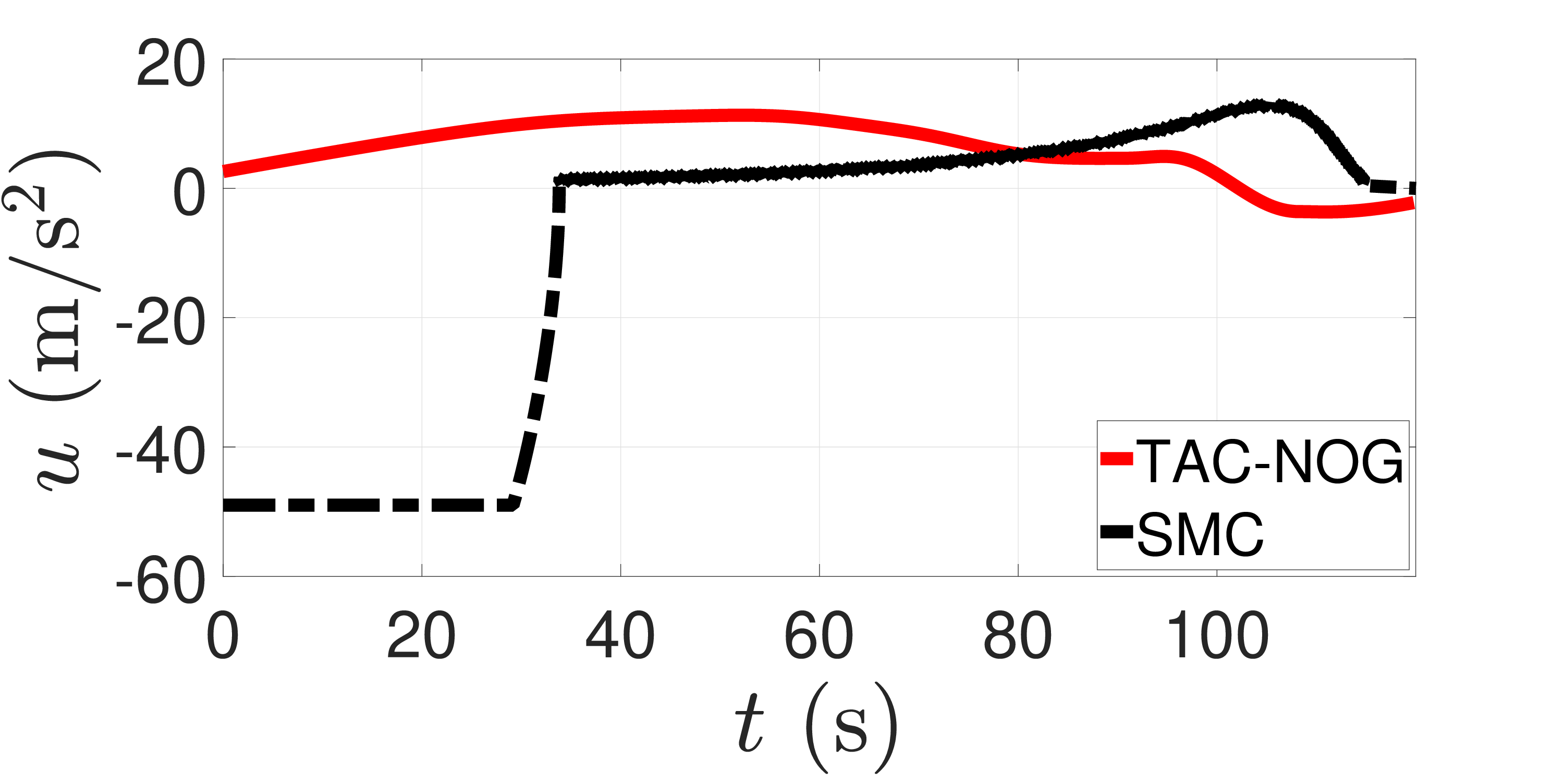} }
        \hspace{-3mm}
    \subfigure[Pursuer $\#4$]{
        \centering
        \includegraphics[width=0.22\textwidth]{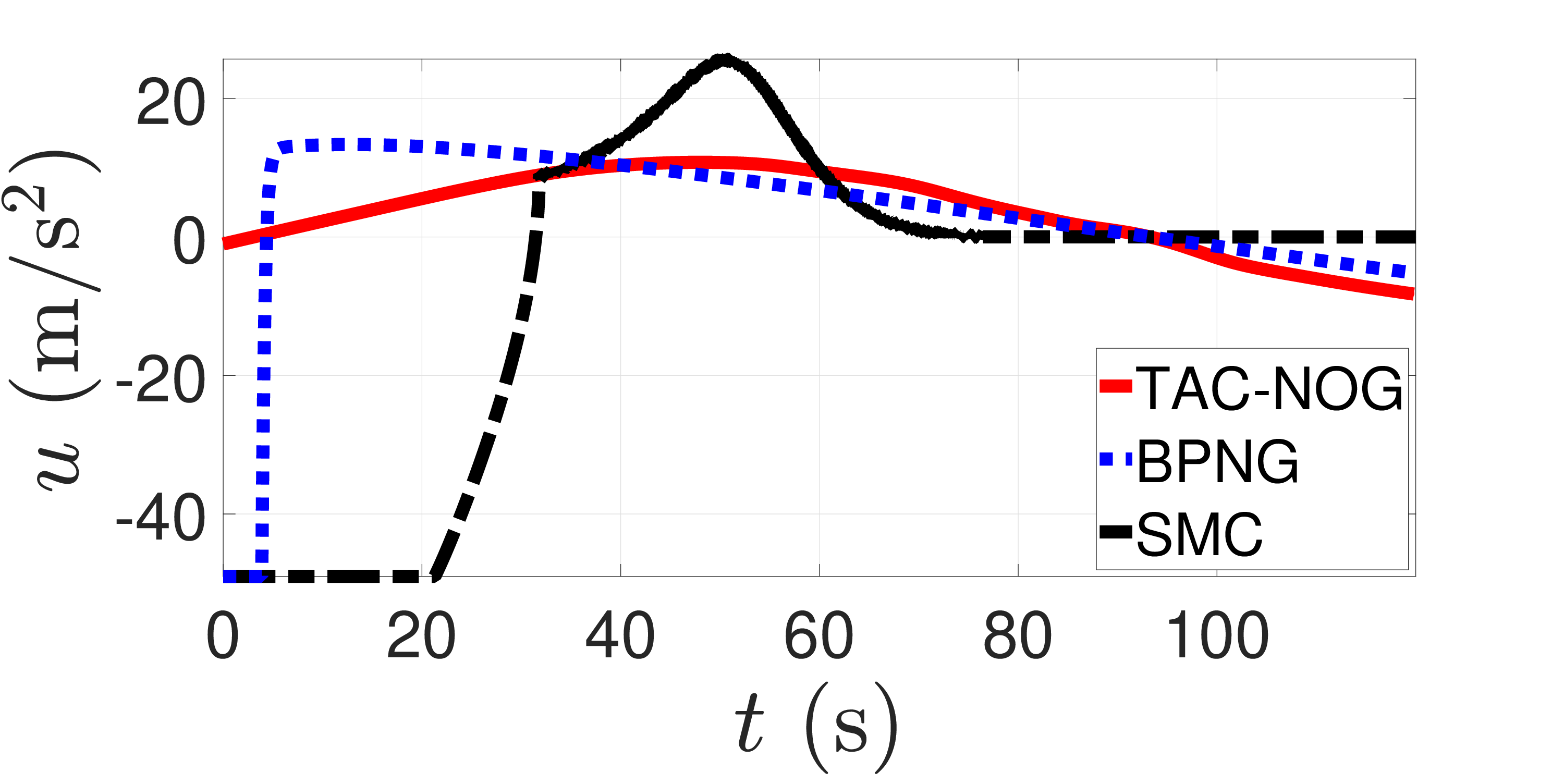}            
    }
    \caption{Case B: control profiles related to three guidance laws.}
    \label{Fig:different angle2}
\end{figure}

Meanwhile, the control efforts required by three methods are quite different. The corresponding control profiles are presented in Fig.~\ref{Fig:different angle2}, and the control effort for each pursuer is listed in Table~\ref{Tab:table2}, where $J_{SMC}$ denotes the control effort required by the SMC-based law and the term “NA” indicates that the corresponding guidance law is not applicable. It is clear from Table~\ref{Tab:table2} that both the BPNG law and the SMC-based guidance law necessitate greater control effort for each pursuer compared to the TAC-NOG law. In order to diminish the control effort of the SMC-based guidance law, one needs to design an appropriate sliding surface. However, dynamically selecting an appropriate sliding surface on board is well-known to be challenging \cite{2015Tokat}. In comparison, the TAC-NOG law developed in the paper can generate the optimal guidance command with expected impact time and impact angle without any tailored operations.

\begin{table}[!hbt]
\centering
\caption{Case B: the control effort required by three guidance laws for different impact angles.}
\label{Tab:table2}
\tiny
\begin{tabular}{c|ccc}
\hline
& $J_{TAC-NOG}~ \rm{(m^2/s^3)}$ & $J_{BPNG}~\rm{(m^2/s^3)}$ & $J_{SMC}~\rm{(m^2/s^3)}$ \\\hline
$P_1$ & $6.265\times10^3$ & $7.685\times10^3$ & $2.570\times10^4$ \\
$P_2$ & $1.791\times10^3$ & NA & NA \\
$P_3$ & $3.526\times10^3$ & NA & $3.941\times10^4$ \\
$P_4$ & $2.808\times10^3$ & $8.671\times10^3$ & $3.545\times10^4$ \\
\hline
\end{tabular}
\end{table}

\section{Conclusions}
\label{Se:Conclusions}

The nonlinear optimal guidance problem with constraints on impact time and impact angle was addressed in this paper by combining the usual optimal control method with a neural network. The PMP was first used to convert the solution trajectory of the minimum-effort control problem into a solution of a parameterized system. Then, by analyzing the geometric property for the solution of the parameterized system, a new optimality condition was established (cf. Lemma \ref{Le:3}). By embedding this optimality condition and the usual dis-conjugacy condition into the parameterized system, a large number of optimal trajectories (at least locally optimal) can be generated by solving some initial value problems. Furthermore, a scaling invariance property was found for the solutions of the parameterized system (cf. Lemma \ref{Le:zoom}). As a consequence,  a simple neural network trained by the solutions of the parameterized system at any fixed time can be used to generate the TAC-NOG command within milliseconds. In addition to the real-time property, the control efforts required by the TAC-NOG command are relatively smaller than those required by optimization methods and existing guidance laws, as demonstrated by the numerical examples.

\begin{ack}                            
The research was supported by the National Natural Science Foundation of China under Grant No.62088101.

\end{ack}

\bibliographystyle{unsrt}        
\bibliography{autosam}           

\begin{thebibliography}{10}

\bibitem{2006Insoo}
In-Soo Jeon, Jin-Ik Lee, and Min-Jea Tahk.
\newblock Impact-time-control guidance law for anti-ship missiles.
\newblock {\em IEEE Transactions on Control Systems Technology}, 14(2):260--266, 2006.

\bibitem{2004Kim}
Sangkeun Jeong, Sungjin Cho, and E.G. Kim.
\newblock Angle constraint biased png.
\newblock In {\em 2004 5th Asian Control Conference (IEEE Cat. No.04EX904)}, volume~3, pages 1849--1854 Vol.3, 2004.

\bibitem{2007Insoo}
Jin-Ik Lee, In-Soo Jeon, and Min-Jea Tahk.
\newblock Guidance law to control impact time and angle.
\newblock {\em IEEE Transactions on Aerospace and Electronic Systems}, 43(1):301--310, 2007.

\bibitem{2013Zhang}
Youan Zhang, Guoxin Ma, and Aili Liu.
\newblock Guidance law with impact time and impact angle constraints.
\newblock {\em Chinese Journal of Aeronautics}, 26(4):960--966, 2013.

\bibitem{2022Chenye}
Ye~Chen, Shufan Wu, and Xiaoliang Wang.
\newblock Impact time and angle control optimal guidance with field-of-view constraint.
\newblock {\em Journal of Guidance, Control, and Dynamics}, 45(12):2369--2378, 2022.

\bibitem{2019ChenZ}
Zheng Chen and Tal Shima.
\newblock Nonlinear optimal guidance for intercepting a stationary target.
\newblock {\em Journal of Guidance, Control, and Dynamics}, 42(11):2418--2431, 2019.

\bibitem{2016Saleem}
Abdul Saleem and Ashwini Ratnoo.
\newblock Lyapunov-based guidance law for impact time control and simultaneous arrival.
\newblock {\em Journal of Guidance, Control, and Dynamics}, 39(1):164--173, 2016.

\bibitem{KIM20142509}
Mingu Kim and Youdan Kim.
\newblock Lyapunov-based pursuit guidance law with impact angle constraint.
\newblock {\em IFAC Proceedings Volumes}, 47(3):2509--2514, 2014.
\newblock 19th IFAC World Congress.

\bibitem{2012Harl}
Nathan Harl and S.~N. Balakrishnan.
\newblock Impact time and angle guidance with sliding mode control.
\newblock {\em IEEE Transactions on Control Systems Technology}, 20(6):1436--1449, 2012.

\bibitem{2020Han}
Tuo Han, Yong Xi, Guangshan Chen, and Qinglei Hu.
\newblock Three-dimensional impact time and angle guidance via controlling line-of-sight dynamics.
\newblock In {\em 2020 Chinese Control And Decision Conference (CCDC)}, pages 2850--2855, 2020.

\bibitem{2016Zhao}
Yao Zhao, Yongzhi Sheng, and Xiangdong Liu.
\newblock Analytical impact time and angle guidance via time-varying sliding mode technique.
\newblock {\em ISA Transactions}, 62:164--176, 2016.
\newblock SI: Control of Renewable Energy Systems.

\bibitem{2019Chenxiaot}
Xiaotian Chen and Jinzhi Wang.
\newblock Sliding-mode guidance for simultaneous control of impact time and angle.
\newblock {\em Journal of Guidance, Control, and Dynamics}, 42(2):394--401, 2019.

\bibitem{2006Vadim}
Vadim Utkin and Hoon Lee.
\newblock Chattering problem in sliding mode control systems.
\newblock {\em IFAC Proceedings Volumes}, 39(5):1, 2006.
\newblock 2nd IFAC Conference on Analysis and Design of Hybrid Systems.

\bibitem{2013Kim}
Tae-Hun Kim, Chang-Hun Lee, In-Soo Jeon, and Min-Jea Tahk.
\newblock Augmented polynomial guidance with impact time and angle constraints.
\newblock {\em IEEE Transactions on Aerospace and Electronic Systems}, 49(4):2806--2817, 2013.

\bibitem{2016Tekin}
Koray~S. Erer and Raziye Tekin.
\newblock Impact time and angle control based on constrained optimal solutions.
\newblock {\em Journal of Guidance, Control, and Dynamics}, 39(10):2448--2454, 2016.

\bibitem{2019Kang}
Shen Kang, Raziye Tekin, and Florian Holzapfel.
\newblock Generalized impact time and angle control via look-angle shaping.
\newblock {\em Journal of Guidance, Control, and Dynamics}, 42(3):695--702, 2019.

\bibitem{2012Harrison}
Gregg~A. Harrison.
\newblock Hybrid guidance law for approach angle and time-of-arrival control.
\newblock {\em Journal of Guidance, Control, and Dynamics}, 35(4):1104--1114, 2012.

\bibitem{2015Kumar2}
Shashi~Ranjan Kumar and Debasish Ghose.
\newblock Impact time and angle control guidance.
\newblock In {\em AIAA guidance, navigation, and control conference}, page 0616, 2015.

\bibitem{2017Song}
Junhong Song, Shenmin Song, and Shengli Xu.
\newblock Three-dimensional cooperative guidance law for multiple missiles with finite-time convergence.
\newblock {\em Aerospace Science and Technology}, 67:193--205, 2017.

\bibitem{2018Hu}
Qinglei Hu, Tuo Han, and Ming Xin.
\newblock New impact time and angle guidance strategy via virtual target approach.
\newblock {\em Journal of Guidance, Control, and Dynamics}, 41(8):1755--1765, 2018.

\bibitem{1998Shneydor}
N~A Shneydor.
\newblock {\em Missile Guidance and Pursuit}.
\newblock Missile Guidance and Pursuit, 1998.

\bibitem{2022Merkulov}
Gleb Merkulov, Martin Weiss, and Tal Shima.
\newblock Minimum-effort impact-time control guidance using quadratic kinematics approximation.
\newblock {\em Journal of Guidance, Control, and Dynamics}, 45(2):348--361, 2022.

\bibitem{1984Guelman}
M.~Guelman and J.~Shinar.
\newblock Optimal guidance law in the plane.
\newblock {\em Journal of Guidance, Control, and Dynamics}, 7(4):471--476, 1984.

\bibitem{2022Wang}
Kun Wang, Zheng Chen, Han Wang, Jun Li, and Xueming Shao.
\newblock Nonlinear optimal guidance for intercepting stationary targets with impact-time constraints.
\newblock {\em Journal of Guidance, Control, and Dynamics}, 45(9):1614--1626, 2022.

\bibitem{Cheng2023}
Lin Cheng, Han Wang, Shengping Gong, and Xu~Huang.
\newblock Neural-network-based nonlinear optimal terminal guidance with impact angle constraints.
\newblock {\em IEEE Transactions on Aerospace and Electronic Systems}, pages 1--11, 2023.

\bibitem{2006Luping}
Ping Lu and Frank Chavez.
\newblock Nonlinear optimal guidance.
\newblock In {\em AIAA Guidance, Navigation, and Control Conference and Exhibit}, page 6079, 2006.

\bibitem{1987Pontryagin}
Lev~Semenovich Pontryagin.
\newblock {\em Mathematical theory of optimal processes}.
\newblock Routledge, 2018.

\bibitem{1989Kurt}
Kurt Hornik, Maxwell Stinchcombe, and Halbert White.
\newblock Multilayer feedforward networks are universal approximators.
\newblock {\em Neural Networks}, 2(5):359--366, 1989.

\bibitem{2018Izzo2}
Carlos S\'{a}nchez-S\'{a}nchez and Dario Izzo.
\newblock Real-time optimal control via deep neural networks: Study on landing problems.
\newblock {\em Journal of Guidance, Control, and Dynamics}, 41(5):1122--1135, 2018.

\bibitem{2016Sanchez}
Carlos Sánchez-Sánchez, Dario Izzo, and Daniel Hennes.
\newblock Learning the optimal state-feedback using deep networks.
\newblock In {\em 2016 IEEE Symposium Series on Computational Intelligence (SSCI)}, pages 1--8, 2016.

\bibitem{2016chen}
Z.~Chen, J.-B. Caillau, and Y.~Chitour.
\newblock $\mathrm{L^1}$-minimization for mechanical systems.
\newblock {\em SIAM Journal on Control and Optimization}, 54(3):1245--1265, 2016.

\bibitem{2015Tokat}
Sezai Tokat, M.~Sami Fadali, and Osman Eray.
\newblock {\em A Classification and Overview of Sliding Mode Controller Sliding Surface Design Methods}, pages 417--439.
\newblock Springer International Publishing, Cham, 2015.

\end{thebibliography}



\appendix

\section{The rotation of trajectory}
\label{App:A}

Given a trajectory $(x(t),y(t),\theta(t))$ of Eq.~(\ref{Eq:system1}) for $t\in [0,t_f]$ with $(x(t_f),y(t_f)) = (0,0)$, let $u(t)$ for $t\in [0,t_f]$ be the corresponding control.
Then, for any $\psi \in [0,2\pi]$, there exists a trajectory $(\bar{x}(t),\bar{y}(t),\bar{\theta(t)})$ for $t\in [0,t_f]$ so that
\begin{equation}
\label{Eq:rotation1}
    \begin{aligned}
        [\bar{x}(t),\bar{y}(t),\bar{\theta}(t)]^T = &
        \left[
        \begin{array}{ccc}
        \cos \psi & -\sin \psi & 0\\
        \sin \psi & \cos \psi & 0\\
        0   & 0 &   1 
        \end{array}
        \right][x(t),y(t),\theta(t)]^T\\
        &+[0,0,\psi]^T, \ t\in [0,t_f]
    \end{aligned}
\end{equation}
Rewrite Eq.~(\ref{Eq:rotation1}) as
\begin{equation}
    \label{Eq:rotation2}
    \left\{
        \begin{aligned}
		\bar{x}(t) &= x(t)\cos\psi-y(t)\sin\psi\\
		\bar{y}(t) &= x(t)\cos\psi+y(t)\cos\psi\\
		\bar{\theta}(t)&=\theta(t)+\psi
		\end{aligned}
	\right.
\end{equation}
By differentiating Eq.~(\ref{Eq:rotation2}), the kinematics of the pursuer can be expressed as
\begin{align}
    &\begin{aligned}\label{Eq:rotation3_1}
    	\dot{\bar{x}} &= \dot{x}\cos\psi-\dot{y}\sin\psi+(-x(t)\sin\psi-y(t)\cos\psi)\dot{\psi}\\
     &=\cos\theta(t)\cos\psi-\sin\theta(t)\sin\psi\\
     &=\cos \bar{\theta}(t)
\end{aligned}\\
 &\begin{aligned}\label{Eq:rotation3_2}
    	\dot{\bar{y}} &= \dot{x}\sin\psi+ \dot{y}\cos\psi+(x(t)\cos\psi-y(t)\sin\psi)\dot{\psi}\\
     &=\cos\theta(t)\sin\psi+\sin\theta(t)\cos\psi\\
     &=\sin \bar{\theta}(t)
\end{aligned}\\
&\begin{aligned}\label{Eq:rotation3_3}
    	\dot{\bar{\theta}} &= \dot{\theta}+\dot{\psi}=u(t)
\end{aligned}
\end{align}
where, for any constant $\psi\in[0,2\pi)$, $\dot{\psi}=0$ always holds. Upon comparing Eqs.~(\ref{Eq:rotation3_1}-\ref{Eq:rotation3_3}) with Eq.~(\ref{Eq:system1}), it is evident that the trajectory $(\bar{x}(t),\bar{y}(t),\bar{\theta(t)})$ for $t\in [0,t_f]$ serves as a solution of Eq.~(\ref{Eq:system1}).
Therefore, we can use the coordinate rotation in Eq.~(\ref{Eq:rotation1}) to change the final heading angle while keeping the control $u$ unchanged.

\section{Proofs of Lemmas in Section \ref{Se:Characterizations}}
\label{App:C}

\begin{figure}
    \centering
    \subfigure[]{
        \centering
        \includegraphics[height=4cm]{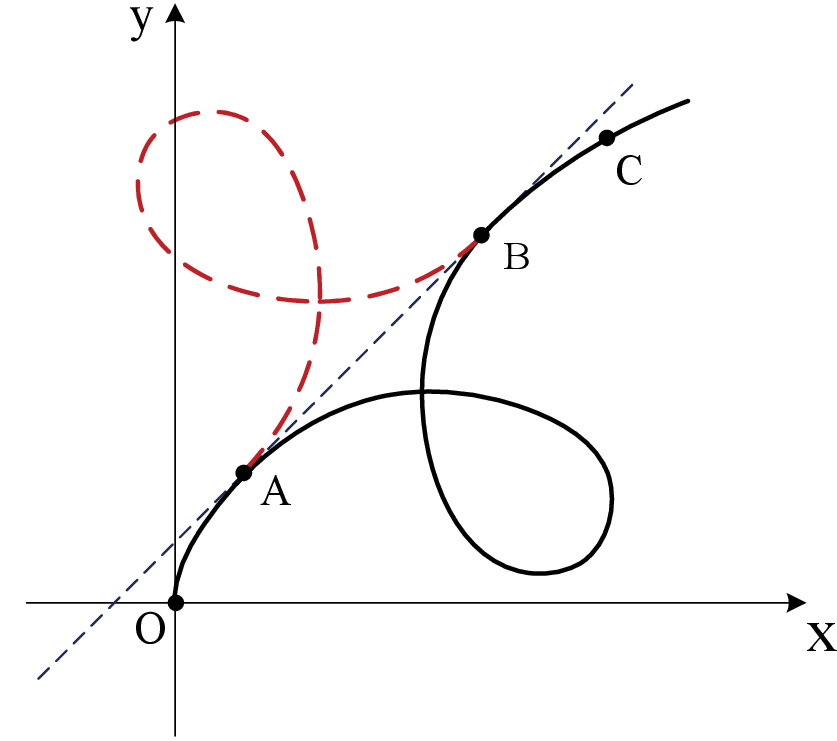}
   }
    ~~~~
    \subfigure[]{
        \centering
        \includegraphics[height=4cm]{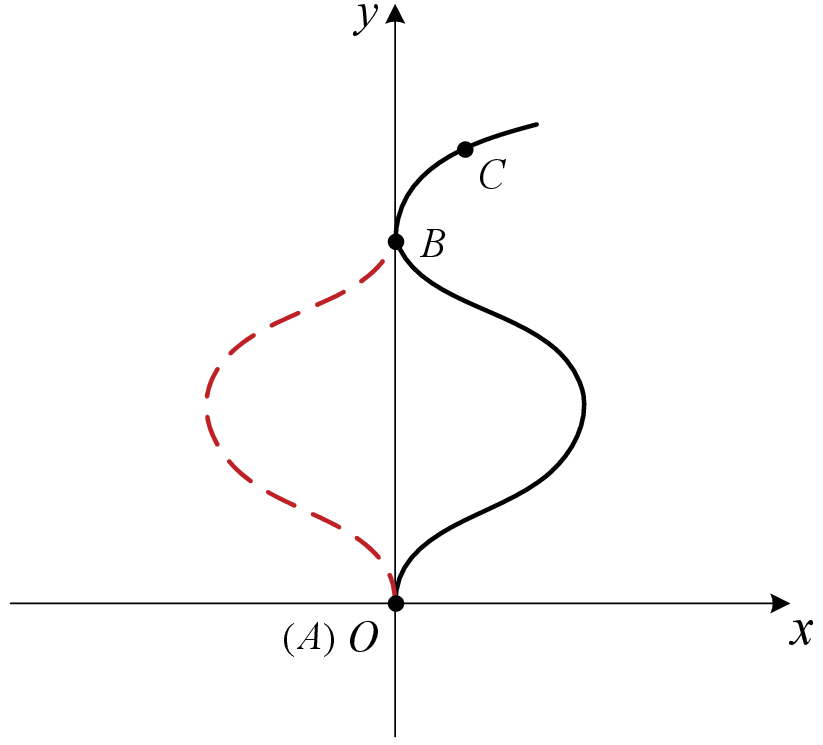}
    }
    \caption{Geometry for the co-linear case.}
    \label{Fig:a_1}
\end{figure}

Proof of Lemma \ref{Le:3}. By contradiction, assume that there are two states A and B along an extremal trajectory $\boldsymbol{Z}(t,\bar{\boldsymbol{q}})$ for $t\in[0,t_f]$, where the velocity vectors of these two states are co-linear, as shown in Fig.~\ref{Fig:a_1}. Denote by $t_1\in$ $(0,t_f)$ and $t_2\in$ $(0,t_1)$ the time when the pursuer reaches state A and state B, respectively, i.e., $A=\boldsymbol{Z}(t_1,\bar{\boldsymbol{q}})$ and $B=\boldsymbol{Z}(t_2,\bar{\boldsymbol{q}})$.
Let $u(t)$ for $t\in[0,t_f]$ be the corresponding control command along the trajectory. 

Set $\overline{u}=-u(t)$ for $t\in [t_2,t_1]$. Then, the trajectory of the pursuer steered by control $\overline{u}$ on $[t_2,t_1]$ is symmetric to that steered by $u$ with respect to the common tangent of state A and state B, as shown in Fig.~\ref{Fig:a_1}, where the solid curves denote the extremal trajectory $\boldsymbol{Z}(t,\bar{\boldsymbol{q}})$ on $[0,t_f]$, and the dashed curves denote the trajectory associated with $\overline{u}=-u(t)$ for $t\in[t_2,t_1]$.

Let us choose a time $\overline{t}\in [0,t_2)$, and let C denote the state at $\overline{t}$. 
Denote by $\gamma$ the piece of extremal trajectory $\boldsymbol{Z}(t,\bar{\boldsymbol{q}})$ from C to the origin, and denote by $\overline{\gamma}$ the smooth concatenation of extremal trajectory $\boldsymbol{Z}(t,\bar{\boldsymbol{q}})$ from C to B, the extremal trajectory of the dashed curve, and the extremal trajectory $\boldsymbol{Z}(t,\bar{\boldsymbol{q}})$ from A to origin.
Denote by $t_2^-$ (resp. $t_2^+$) the time that $t \to t_2$ from the direction that $t<t_2$ (resp. $ t>t_2$) when moving along the trajectory $\overline{\gamma}$. According to the definition of trajectory $\overline{\gamma}$, the sign of the control $u(t_2^-)$ is opposite to that of the control $u(t_2^+)$. Therefore, the control command is discontinuous at B when moving along the trajectory $\overline{\gamma}$.
This contradicts the necessary condition in Eq.~(\ref{Eq:ut}) that the control is continuous.
Therefore, the trajectory $\overline{\gamma}$ is not optimal. However, it is obvious that the cost for the pursuer to move along trajectory $\gamma$ is the same as the cost along trajectory $\overline{\gamma}$. That is to say, trajectory $\gamma$ is not optimal as well.
Therefore, there is another extremal trajectory from point C to the origin with a duration of $t_f - \overline{t}$ and a terminal impact angle $\theta_f=-\pi/2$, resulting in a smaller cost. This completes the proof. 

Proof of Lemma \ref{Le:zoom}. 
Since $(x(t),y(t),\theta(t))$ for $t\in[0,t_g]$ represents an optimal trajectory of the pursuer, its kinematics are expressed by Eq.~(\ref{Eq:system1}).
For any $T\in (0,t_g)$, set $\tau=t\dfrac{T}{t_g}$. Then, we have 
\begin{equation}
    \label{Eq:systemB2}
    \left\{
        \begin{aligned}
		x(\tau)&=x(\tau \dfrac{t_g}{T}) \dfrac{T}{t_g}\\
        y(\tau)&=y(\tau \dfrac{t_g}{T}) \dfrac{T}{t_g}\\
		\theta(\tau)&=\theta(\tau \dfrac{t_g}{T})
		\end{aligned}
	\right.
\end{equation}
Substitute Eq.~(\ref{Eq:systemB2}) into Eq.~(\ref{Eq:system1}), we have
\begin{equation}
    \label{Eq:systemB3}
    \left\{
        \begin{aligned}
        \dfrac{dx(\tau)}{d\tau}&=\dfrac{dx(t)}{dt} \dfrac{T}{t_g} \dfrac{dt}{d\tau} =\cos\theta(t)\\
        \dfrac{dy(\tau)}{d\tau}&=\dfrac{dy(t)}{dt} \dfrac{T}{t_g} \dfrac{dt}{d\tau} =\sin\theta(t)\\
        \dfrac{d\theta(\tau)}{d\tau} &= \dfrac{d\theta(t)}{dt} \dfrac{dt}{d\tau} =u(t)\dfrac{t_g}{T}
		\end{aligned}
	\right.
\end{equation}
Thus, the following relationship can be obtained
\begin{align}
	u(\tau)=u(t)\dfrac{t_g}{T}
\end{align}
When $\tau=0$, the optimal feedback control associated with the current state $(x_c,y_c,\theta_c)$ can be expressed as 
\begin{align}
\label{Eq:ut2}
	u^*(t_g,x_c,y_c,\theta_c)=\dfrac{T}{t_g} u^*(T,\dfrac{T}{t_g} x_c,\dfrac{T}{t_g} y_c,\theta_c)
\end{align}
which completes the proof.

Proof of Lemma \ref{Le:chen}.
According to Theorem 2 in \cite{2019ChenZ}, for any $\boldsymbol{z}\in \mathcal{Z}_T$, there exists an optimal trajectory $(x(t),y(t),\theta(t))$ for $t\in [0,T]$ of Problem \ref{Problem1} so that $\boldsymbol{z} = (x(0),y(0),\theta(0))$. In addition, by the definition of the $\boldsymbol{q}$-parameterized system in Eq.~(\ref{Eq:system2}), there exists $\boldsymbol{q}$ so that $(x(t),y(t),\theta(t))= Z(t,\boldsymbol{q})$ for $t\in [0,T]$, which completes the proof of the first statement. 

Then, let us prove the second statement. For any given $\boldsymbol{q}$, an extremal trajectory $\boldsymbol{Z}(t,\boldsymbol{q})$ on $[0,t_f]$ with a final condition of $(0,0,-\pi/2)$ can be obtained by integrating the $\boldsymbol{q}$-parameterized system in Eq.~(\ref{Eq:system2}).
By the definition of $\mathcal{Z}_T$, it is clear that by setting the value of integral duration as $T$, $\boldsymbol{Z}(T,\boldsymbol{q})$ is exactly the state in set $\mathcal{Z}_T$.
Thus, for any $\boldsymbol{q}\in \mathbb{R}^3$, $Z(T,\boldsymbol{q})\in \mathcal{Z}_T$ always holds, completing the proof of the second statement of Lemma \ref{Le:chen}.

\end{document}